\documentclass[10pt]{amsart}

\usepackage{amssymb}
\usepackage{verbatim,setspace,graphicx}
\usepackage[colorlinks=true, citecolor=blue]{hyperref}
\usepackage{multirow,array}

\newtheorem{theorem}{Theorem}

\newtheorem{proposition}[theorem]{Proposition}

\theoremstyle{remark}
\newtheorem{remark}{Remark}[theorem]

\newcommand{\T}		{\mathbb{T}}
\newcommand{\D}		{\mathbb{D}}

\newcommand{\C}		{\mathbb{C}}
\newcommand{\N}		{\mathbb{N}}
\newcommand{\Z}		{\mathbb{Z}}
\newcommand{\Om}		{\mathbb{O}}
\newcommand{\map}	{\Phi}

\newcommand{\const}	{\textnormal{const.}}

\newcommand{\dist}	{\textnormal{dist}}

\newcommand{\re}	{\textnormal{Re}}

\begin{document}

\title[Universality for potential theoretic ensembles]{Universality for ensembles of matrices with potential
  theoretic weights on domains with smooth boundary}

\author[C. Sinclair]{Christopher D. Sinclair}

\address{Department of Mathematics, University of Oregon, Eugene, OR, 97403}

\email{csinclai@uoregon.edu}

\author[M. Yattselev]{Maxim L. Yattselev}

\address{Corresponding author, Department of Mathematics, University of Oregon, Eugene, OR, 97403}

\email{maximy@uoregon.edu}

\subjclass[2000]{}

\keywords{normal matrix model, universality, scaling limits, exterior
  asymptotics, area orthogonality}

\maketitle

\begin{abstract}
  We investigate a two-dimensional statistical model of $N$ charged
  particles interacting via logarithmic repulsion in the presence of
  an oppositely charged compact region $K$ whose charge density is
  determined by its equilibrium potential at an inverse temperature
  corresponding to $\beta = 2$. When the charge on the region, $s$, is
  greater than $N$, the particles accumulate in a neighborhood of the
  boundary of $K$, and form a determinantal point process on
  the complex plane.  We investigate the scaling limit, as
  $N \rightarrow \infty$, of the associated kernel in the neighborhood
  of a point on the boundary under the assumption that the boundary is
  sufficiently smooth.  We find that the limiting kernel depends on the
  limiting value of $N/s$, and prove universality for these kernels.
  That is, we show that, the scaled kernel in a neighborhood of a
  point $\zeta \in \partial K$ can be succinctly expressed in terms of
  the scaled kernel for the closed unit disk, and the
  exterior conformal map which carries the complement of $K$ to the
  complement of the closed unit disk.  When $N / s \rightarrow 0$ we
  recover the universal kernel discovered by Lubinsky in \cite{Lub10}.
\end{abstract}

\section{Introduction} 

\subsection{Potential Theoretic Weights} 

Let $K \subseteq \C$ be a compact subset whose boundary
$T = \partial K$ is a Jordan curve.  We will assume that $T$ is
sufficiently nice in a way that will be made precise in the sequel.
For such $K$, there exists a unique measure $\omega_K$, the
\emph{equilibrium measure} on $K$, that minimizes the energy
functional $I[\sigma]:=-\iint\log|z-u|d\sigma(u)d\sigma(z)$ among all
positive probability measures $\sigma$ supported on $K$
\cite{Ransford}. We define $P_K: \C \rightarrow (0,\infty)$, by
\begin{equation}
\label{PK}
P_K(z) := \exp\left\{I[\omega_K]+\int\log|z-u|d\omega_K(u)\right\}.
\end{equation}
This function is simply the rescaled exponentiated {\em equilibrium
  potential} of $K$.
 
It can be verified that $\omega_K$ is supported on $T$, $P_K$ is
identically one on $K$ and, as $z \to \infty$, $P_K(z)/|z|
\to\gamma_K^{-1}$, where $\gamma_K:=\exp\{-I[\omega_K]\}$ is the
\emph{logarithmic capacity} of $K$. 

  In this paper we will be interested in random vectors whose joint
density is given by  
\begin{equation}
\label{eq:1}
\Omega_N(\boldsymbol \lambda) := \frac{1}{Z_N} \bigg\{\prod_{n=1}^N w(\lambda_n)
\bigg\} \prod_{m < n} | \lambda_n - \lambda_m |^2; \qquad \boldsymbol
\lambda 
\in \C^N, 
\end{equation}
where 
\begin{equation}
\label{eq:8}
w(\lambda) = P_K(\lambda)^{-2s}, \qquad 
Z_N := \int_{\C^N} \bigg\{\prod_{n=1}^N w(\lambda_n)
\bigg\} \prod_{m < n} | \lambda_n - \lambda_m |^2 \,
dA^N(\boldsymbol\lambda),
\end{equation}
and $s > N$ (that is $s$ is sufficiently large to guarantee that $Z_N$
is finite).  Here and throughout, $A$ and $A^N$ are Lebesgue measure
on $\C$ and $\C^N$ respectively.  

We will often refer to the components of such random vectors as eigenvalues,
since the joint density (\ref{eq:1}) can be thought of as a
modification of the joint eigenvalue density of the ensemble of
matrices with i.i.d. complex Gaussian entries.  The eigenvalues of
this latter matrix ensemble, 
originally introduced by Ginibre \cite{MR0173726}, have joint density
given by (\ref{eq:1}) where $w(\lambda) = e^{-|\lambda|^2}$.  

In Section~\ref{sec:two-other-models} we will give (i) a matrix model
whose joint density of eigenvalues is given by (\ref{eq:1}) with
weight given as in (\ref{eq:8}) for $K$ equal to the closed unit disk,
as well as (ii) models for more 
general $K$  where the components of $\boldsymbol\lambda$ represent
the positions of electrostatic particles confined to the
plane and in the presence of a field determined by $K$, and (iii)
an ensemble of random polynomials chosen with respect to a height
function determined by $K$ whose roots are distributed as in
(\ref{eq:1}).  

Our primary goal is to demonstrate that, in the double scaling limit
as $s$ and $N$ approach infinity, the local statistics of the
eigenvalues near a point on the boundary of $K$ depend only on the
limiting ratio of $s$ and $N$, but are essentially independent of the
specifics of $K$.  This will follow from the asymptotic behavior of
the reproducing kernel of $L^2(w)$, which in turn follows from the
asymptotics of the leading coefficient of the related orthonormal
polynomials. When $s = \infty$, $w$ is simply the characteristic
function of $K$ and our results collapse to those given by Lubinsky
\cite{Lub10} for the universality of reproducing kernels formed with
respect to Bergman polynomials for $K$.

\subsection{Eigenvalue Statistics}

We briefly review some basic concepts for solvable ensembles of
random matrices and how they relate to eigenvalue statistics.  In this
section we will assume that the joint density of eigenvalues is given
by (\ref{eq:1}) where, for the purposes of this section, $w: \C
\rightarrow [0,\infty)$ is any non-negative function such that $0 <
Z_N < \infty$.  

We will suppose that $\Xi = \{ \xi_1, \xi_2, \ldots, \xi_N \} \subset
\C$ is a random set corresponding to the eigenvalues of a random matrix
from our ensemble.  (Or, what amounts to the same thing, $\Xi$ is the
set corresponding to a random vector sampled from the density
$\Omega_N$).  Given a set $E \subseteq \C$ we may construct a random
variable $X$ given by the cardinality of $\Xi \cap E$.  Given disjoint
subsets $E_1, E_2, \ldots, E_n$ we will let $X_1, X_2, \ldots, X_n$ be
the corresponding random variables.  The $n$-th {\em correlation
  function} of our ensemble is defined to be $R_n: \C^n \rightarrow [0,
\infty)$, where 
\[
E[X_1 X_2 \cdots X_n] := \int_{E_1} \int_{E_2} \cdots \int_{E_n}
R_n(\boldsymbol \lambda) \, dA^n(\boldsymbol \lambda).
\]
It is straightforward to see that $R_N = N! \, \Omega_N$.  A less
obvious exercise is to show that for $0 \leq n \leq
N$, 
\begin{equation}
\label{eq:2}
R_n(\boldsymbol \lambda) = \frac{1}{(N-n)!} \int_{\C^{N-n}}
R_N(\boldsymbol \lambda \vee \mathbf x) \, dA^{N-n}(\mathbf x),
\end{equation}
where $\boldsymbol \lambda \vee \mathbf x = (\lambda_1, \cdots,
\lambda_n, x_1, \ldots, x_{N-n})$.  Many probabilities of interest can
be expressed in terms of correlation functions.  One particularly
important example is the {\em gap} probability that there are no
eigenvalues in $E$,
\begin{equation}
\label{eq:4}
\mathrm{Prob}\{X = 0 \} = \sum_{n=0}^N \frac{(-1)^n}{n!}
\int_{E^n} R_n(\boldsymbol \lambda) \, dA^n(\boldsymbol \lambda).
\end{equation}
Equations (\ref{eq:2}) and (\ref{eq:4}) are valid for a wide variety
of symmetric measures on $\C^N$.  However, the presence of the square
of the Vandermonde determinant which appears in (\ref{eq:1}) leads to
additional structure which may be exploited.

Suppose $\pi_0, \pi_1, \ldots, \pi_{N-1}$ are the orthonormal
polynomials with respect to the weight $w$.  That is,
\[
\int_{\C} \pi_n(z) \overline{\pi_m(z)} w(z) \, dA = \delta_{nm},
\]
where, as usual, $\delta_{nm}$ is 1 or 0 depending on whether or not $n = m$.
The {\em kernel} of the ensemble is defined by 
\[
\widetilde K_N(z, u) := \sqrt{w(z) w(u)} \sum_{n=0}^{N-1} \pi_n(z)
\overline{\pi_n(u)}. 
\]
(Following Lubinsky's notation, we will reserve the symbol $K_N$ for
the unweighted analog of this kernel).  In a celebrated result, Mehta and
Gaudin \cite{Mehta1960420} were able to express the correlation
functions of ensembles with eigenvalue density (\ref{eq:1}) in terms
of determinants of matrices formed from this kernel,
\begin{equation}
\label{eq:3}
R_n(\boldsymbol \lambda) = \det \left[ \widetilde K_N(\lambda_j,
  \lambda_k) \right]_{j,k=1}^n. 
\end{equation}
(See also \cite{TW} for a more modern derivation).

\subsection{Universality}

When $N$ is large we expect that, with high probability, the 
eigenvalues will accumulate in a neighborhood of $\partial K$.
Slightly more precisely, if $\zeta \in \partial K$, then the number of
eigenvalues in a disk of (small) radius $\epsilon$ about $\zeta$ is
proportional to $N$; the constant of proportionality is given by
the integral of the equilibrium measure over the arc of $\partial K$
contained in the disk.  The exact details of this phenomenon will be
explored in a subsequent paper, for now we use this only as intuition
to guess the proper scale on which we expect $\widetilde K_N$ to converge.   

From (\ref{eq:4}) and (\ref{eq:3}), the probability that there
are no eigenvalues in a disk of radius $\epsilon$ centered at $\zeta$
is given by 
\begin{equation}
\label{eq:5}
\sum_{n=0}^N \frac{(-1)^n}{n!}
\int_{\D^n} \det \left[ \epsilon^2 \widetilde K_{N,s}(\zeta + \epsilon
  \lambda_j, \zeta + \epsilon 
  \lambda_k) \right]_{j,k=1}^n \, dA^n(\boldsymbol \lambda).  
\end{equation}
where $\D$ is the disk of radius $1$ centered at the origin.  Here we
have made explicit that the kernel is dependent on $s$ as
well as $N$.  

Under the assumption that there are $\mathcal O(N)$ eigenvalues in a
neighborhood of $\zeta$, then we should scale 
$\epsilon$ like $1/N$ in order for (\ref{eq:5}) to approach a
non-trivial limit.  That is, the limiting gap probability of there
being no eigenvalues in a shrinking neighborhood with radius 
$\epsilon = 1/N$ is given by
\begin{equation}
\label{eq:6}
\lim_{N \rightarrow \infty} \sum_{n=0}^N \frac{(-1)^n}{n!}
\int_{\D^n} \det \left[ \frac{1}{N^2} \widetilde K_{N,s}\bigg(\zeta +
  \frac{\lambda_j}{N}, \zeta + \frac{\lambda_k}{N}\bigg)
\right]_{j,k=1}^n \, dA^n(\boldsymbol \lambda).
\end{equation}
Since $s > N$, this limit also depends on how $s$ scales with $N$ and
we will assume that $N/s$ converges to some $\ell \in [0,1]$.

If it can be shown that there is some limiting kernel
$\widetilde H_{\zeta, \ell}$ so that 
\[
\frac{1}{N^2} \widetilde K_{N,s}\bigg(\zeta +
  \frac{z}{N}, \zeta + \frac{u}{N}\bigg) \rightarrow \widetilde
  H_{\zeta,\ell}(z,u) 
\]
uniformly on compact subsets of $\C \times \C$, then (\ref{eq:6})
converges to 
\[
\sum_{n=0}^{\infty} \frac{(-1)^n}{n!} \int_{\D^n} \det\left[
  \widetilde H_{\zeta, \ell}(\lambda_j, \lambda_k) \right]_{j,k=1}^n \,
dA^n(\boldsymbol \lambda).
\]
(See for instance \cite[\S 3.4]{MR2760897}).  Our primary result here
is that $\widetilde H_{\zeta,\ell}$ exists, and is dependent on
$\zeta$ and $K$ in only the most trivial manner.  More specifically, we will
express $\widetilde H_{\zeta,\ell}$ in terms of the limiting kernel
for the ensemble formed from the closed unit disk and the value of a
conformal map from $\C \setminus K$ to $\C \setminus \D$ evaluated at
$\zeta$.  This is what is called {\em universality} for potential
theoretic ensembles.

We will also demonstrate that $\widetilde H_{\zeta,\ell}$ is a convex
combination of $\widetilde H_{\zeta, 0}$ (Lubinsky's
limiting kernel) and $\widetilde H_{\zeta, 1}$. 

\subsection{Potential Theoretic Orthogonal Polynomials}

We denote the orthonormal polynomials for the weight $P_K^{-2s}$, $s>1$, by
$\{\pi_{n,s}\}_{n=0}^{\lfloor s-2\rfloor}$.  That is, these are polynomials with positive leading coefficients that satisfy
\begin{equation}
\label{orthogonality}
\int_\C\pi_{n,s}(z)\overline{\pi_{m,s}(z)}P_K^{-2s}(z) \, dA = \delta_{nm}.
\end{equation}
The reproducing kernel for this system of polynomials is given by
\begin{equation}
\label{kernel}
K_{N,s}(z,u) := \sum_{n=0}^{N-1}\pi_{n,s} (z) \overline{\pi_{n,s}(u)}, \quad N\leq\lfloor s-1\rfloor,
\end{equation}
with the weighted kernel given by
\begin{equation}
\label{tildekernel}
\widetilde K_{N,s}(z,u) := P_K^{-s}(z)P_K^{-s}(u)K_{N,s}(z,u).
\end{equation}
Our derivation of $\widetilde H_{\zeta,\ell}$ will follow from the
asymptotics of $K_{N,s}$, which in turn will follow from the
asymptotics of the orthogonal polynomials.  These latter asymptotics
are of independent interest, and they provide the other primary
results of the paper.

\section{Statement of Results}

In what follows, we assume that $T:=\partial K$ is a rectifiable
Jordan curve which is either analytic or of class $C^{p+1,\alpha}$,
where $p$ is a nonnegative integer and $\alpha\in(0,1)$.  That is, the
arclength function of $T$ is $p$ times continuously differentiable as
a periodic function on the real line and its $p$-th derivative is
$\alpha$-H\"older continuous. Denote by $\map$ the conformal map of
$O:=\overline\C\setminus K$ onto $\Om:= \overline\C\setminus\D$ such that
$\map(\infty)=\infty$ and $\map^\prime(\infty)>0$. In the case where
$T$ is an analytic Jordan curve we denote by $\rho(T)<1$ a number such
that $\map^{-1}$ has a univalent extension into
$|w|>\rho(T)$. Moreover, we put $O_\rho:=\map^{-1}(\{|w|>\rho\})$ for
each $\rho>\rho(T)$.

\begin{figure}[h!]
  \centering
  \includegraphics[scale=1.1]{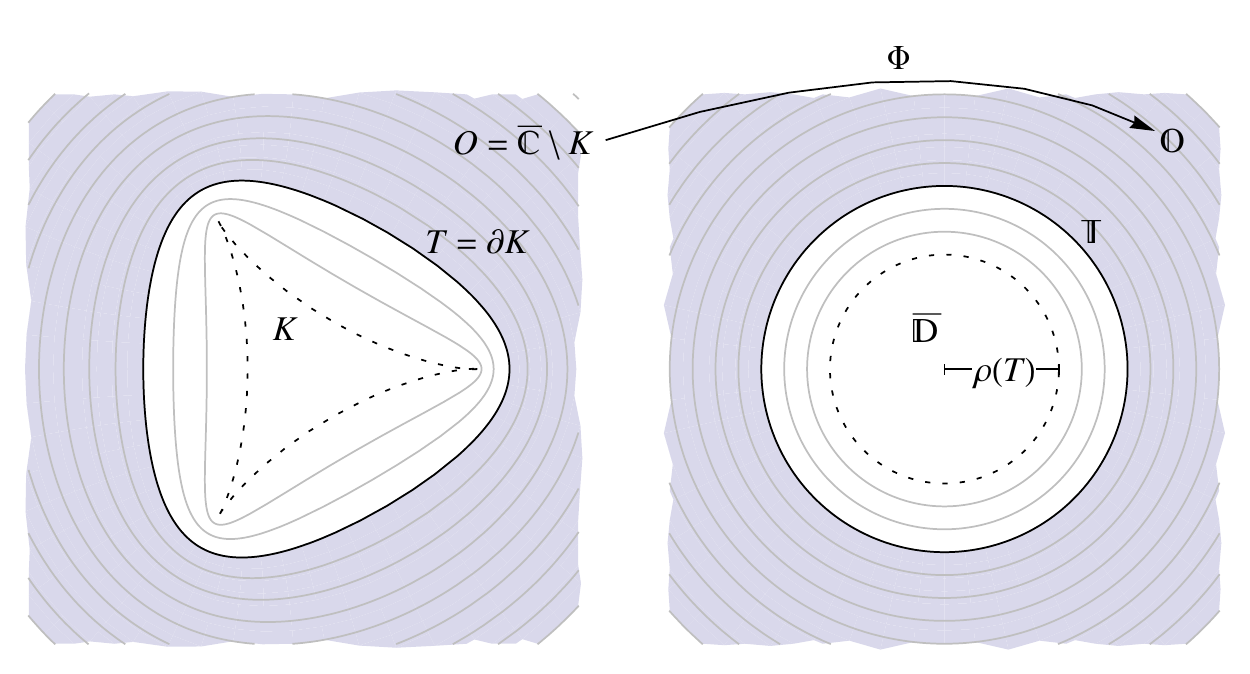}
  \caption{In this figure, $K$ (respectively $\overline \D$) is the
    region enclosed by the black contour.  Here, $K$ has analytic
    boundary, and the dashed contour on the left is the inner-most contour outside
    of which we can find a univalent extension of
    $\Phi^{-1}$. $O_{\rho}$ is represented by the region outside of  
    the contour corresponding to radius $\rho$. The curves outside of
    $T$ are level lines of $P_K$.}
\label{fig:8}
\end{figure}

It is known that $|\map|$ is identically equal to $P_K$ on $O$ and
therefore 
$\map^\prime(\infty)=\gamma_K^{-1}$. Hence, orthogonality relations
\eqref{orthogonality} can be rewritten as
\begin{equation}
\label{anotherform} \int_D\pi_{n,s}(z)\overline{\pi_{m,s}(z)}dA +
\int_O\pi_{n,s}(z)\overline{\pi_{m,s}(z)}|\map(z)|^{-2s}dA =
\delta_{nm},
\end{equation} $n,m\leq\lfloor s-2\rfloor$, where $D$ is the interior domain of
$T$. Since $|\map|>1$ in $O$, we can formally set $\pi_{n,\infty}$ to
be polynomials satisfying
\[ \int_D\pi_{n,\infty}(z)\overline{\pi_{m,\infty}(z)} \, dA =
\delta_{nm}.
\] In a sense, potential theoretic polynomials $\pi_{n,s}$ can be
considered as perturbations of $\pi_{n,\infty}$. The latter were
initially studied by Carleman \cite{Carl23} who derived their exterior
asymptotics (asymptotics in $O$) for the case $T$ being an analytic
Jordan curve. The results in \cite{Carl23} were subsequently extended
by Suetin \cite{Suetin} to include $C^{p+1,\alpha}$ Jordan curves. Other
aspects of the behavior of $\pi_{n,\infty}$, such as zero distribution
and interior asymptotics, were investigated in \cite{Mina08,
DragMina10, DragMina10a}. The following theorem provides an analog of
\cite[Theorem~1.2]{Suetin} for potential theoretic polynomials
$\pi_{n,s}$.

\begin{theorem}
\label{thm:polynomials} 
Let $T=\partial K$ be a Jordan curve of class $C^{p+1,\alpha}$,
$p+\alpha>1/2$, and $\{\pi_{n,s}\}_{n=0}^{\lfloor s-2\rfloor}$ be a
sequence of polynomials satisfying orthogonality relations 
\eqref{orthogonality}. Then, as $n, s \rightarrow \infty$, the leading
coefficient $\varkappa_{n,s}$ of $\pi_{n,s}$ satisfies
\begin{equation}
\label{kappa}
\varkappa_{n,s} =
\frac{1}{\gamma_K^{n+1}}\sqrt{\frac{n+1}{\pi}\left(1-\frac{n+1}{s}\right)}\left[1+\mathcal{O}\left(\frac{1}{n^{2(p+\alpha)}}\right)\right].
\end{equation}
Moreover, if $T$ is an analytic Jordan curve,
then the error terms in \eqref{kappa} can be replaced by
$\mathcal{O}(\rho^{2n})$ for any $\rho(T)<\rho<1$. It also holds that
\begin{equation}
\label{asymp}
\pi_{n,s} =\sqrt{\frac{n+1}{\pi}\left(1-\frac{n+1}{s}\right)} \map^n\map^\prime\big[1+\mathcal{O}\left(\Sigma_n\right)\big]
\end{equation}
uniformly on $\overline O$ as $n,s\to\infty$, where $\Sigma_n$ is
given by Table~\ref{table1}.
\end{theorem}

\begin{table}[h]
\begin{tabular}{|c | m{2cm} | m{2cm} | m{1.5cm} | m{1.5cm} | m{2.5cm} | c}
\cline{1-6}
\multirow{6}{*}{$T$} & \multicolumn{1}{c|}{\multirow{2}{*}{Analytic}} & \multicolumn{4}{c|}{\multirow{2}{*}{$C^{p+1,\alpha}$}} & \\[.4cm]
\cline{2-6}
& \multicolumn{1}{c|}{\multirow{4}{*}{$\rho(T)<\rho<1$}} & \multicolumn{1}{c|}{\multirow{4}{*}{$p\geq2$}} & \multicolumn{2}{c|}{\multirow{2}{*}{$\displaystyle\limsup_{n,s\to\infty}n/s<1$}} & \multicolumn{1}{c|}{\multirow{2}{*}{$\displaystyle\limsup_{n,s\to\infty}n/s=1$}} & \\[.4cm]
\cline{4-6}
&&& \multicolumn{1}{c|}{\multirow{2}{*}{$p=1$}} & \multicolumn{1}{c|}{\multirow{2}{*}{$p=0$}} & \multicolumn{1}{c|}{\multirow{2}{*}{$p=1$}} & \\[.4cm]
\cline{1-6}
\multicolumn{1}{|c|}{$\Sigma_n$} &\center$\rho^n$ & \center$\displaystyle\frac{\log n}{n^{p+\alpha}}$& \center$\displaystyle\frac{\log n}{n^{1+\alpha}}$ & \center$n^{1-2\alpha}$ & \center$n^{-2\alpha}$ & \\[.4cm]
\cline{1-6}
\end{tabular}
\caption{The error term $\Sigma_n$ depending on the smoothness of $T$.}
\label{table1}
\end{table}

\begin{remark}
When $p=0$ and $\limsup_{n,s\to\infty}n/s=1$ the authors were unable
to show that $\Sigma_n\to0$ as $n\to\infty$ (the employed method
yields $\Sigma_n=n^{2(1-\alpha)}$), which is the reason this case is
not included in Table~\ref{table1}. 
\end{remark}

In general, the location of the zeros of $\pi_{n,s}$ depends on $s$ (as well as, obviously, $K$ and $n$).  However, as the following proposition shows, this is not the case for a family of ellipses which interpolate between the unit circle and the interval $[-2,2]$.  

\begin{proposition}
\label{prop:ellipse} 
Let $q \in [0,1)$ and define $\phi(w) := \Phi^{-1}(w) = w+\frac qw$ so
that $K$ is the ellipse bounded by $\phi(\T)$.  Then, for all
$n\leq\lfloor s-2\rfloor$ and all $s$ including $s=\infty$, 
\[
\pi_{n,s} =
\sqrt{\frac{n+1}{\pi}\left(1-\frac{n+1}{s}\right)/\left(1-q^{2n+2}\frac{s-n-1}{s+n+1}\right)}\map^n\map^\prime\left(1-\frac{q^{n+1}}{\map^{2n+2}}\right). 
\] 
That is, the polynomials $\pi_{n,s}$ are the
renormalized Chebysh\"ev polynomials of the second kind for the
interval $\left[-2\sqrt q,2\sqrt q\right]$, where $\pm2\sqrt q$ are
the foci of $T$. 
\end{proposition}

\begin{remark}
In the proposition above all the ellipses have unit logarithmic
capacity ({\it i.e.}, $\gamma_K=1$).
\end{remark} 

\begin{remark}
If $q=0$, then $T=\T$ ($K=\overline\D$) and
\[
\pi_{n,s}(z)=\sqrt{\frac{n+1}{\pi}\left(1-\frac{n+1}{s}\right)}z^n.
\]
\end{remark}

\begin{remark}
If $q = 1$, then $T=K=[-2,2]$ and 
\begin{equation}
\label{-22}
\pi_{n,s}(z)=\sqrt{\frac{s^2-(n+1)^2}{2\pi s}}\frac{1}{\sqrt{z^2-4}}\left[\left(\frac{z+\sqrt{z^2-4}}{2}\right)^{n+1}-\left(\frac{z-\sqrt{z^2-4}}{2}\right)^{n+1}\right].
\end{equation}
It is easy to see that the asymptotic behavior of the normalizing constant in \eqref{-22} is different from the one in \eqref{kappa}. However, the case $K=[-2,2]$ is not covered by Theorem~\ref{thm:polynomials}.
\end{remark}

Theorem~\ref{thm:polynomials} is the essential building block in proving results on asymptotic behavior of kernels $K_{N,s}$ and $\widetilde K_{N,s}$ defined in \eqref{kernel} and \eqref{tildekernel}, respectively. 

\begin{proposition}
\label{prop:kernels}
Let $N\leq\lfloor s-1\rfloor$. Under the conditions of Theorem~\ref{thm:polynomials}, it holds that
\begin{eqnarray}
K_{N,s}(z,w) &=& \frac{\map^\prime(z)\overline{\map^\prime(w)}}{\pi}\left[\left(1-\frac{N+1}{s}\right)\left(-(N+1)\frac{\left[\map(z)\overline{\map(w)}\right]^N}{1-\map(z)\overline{\map(w)}}+\frac{1-\left[\map(z)\overline{\map(w)}\right]^{N+1}}{\left[1-\map(z)\overline{\map(w)}\right]^2}\right)\right. \nonumber \\
&& +\left. \frac1s\left((N+2)\frac{1+\left[\map(z)\overline{\map(w)}\right]^{N+1}}{\left[1-\map(z)\overline{\map(w)}\right]^2} - 2\frac{1-\left[\map(z)\overline{\map(w)}\right]^{N+2}}{\left[1-\map(z)\overline{\map(w)}\right]^3}\right) \right] \nonumber \\
\label{kernel2}
&& + \mathcal{O}\left(\max\left\{1,N^2\Sigma_N\right\}\right)
\end{eqnarray}
uniformly for $z,w\in\overline O$, $z\neq w$, $\dist(z,\partial K)\leq\const/N$ and $\dist(w,\partial K)\leq\const/N$. Moreover, it holds that
\begin{equation}
\label{kernel3}
 K_{N,s}(z,z) = \frac{|\map^\prime(z)|^2}{\pi}\left[\frac{N(N+1)}{2}\left(1-\frac{N+1}{s}\right)+\frac{N(N+1)(N+2)}{6s}\right] + \mathcal{O}\left(\max\left\{1,N^2\Sigma_N\right\}\right)
\end{equation}
uniformly for $z\in\partial K$. 
\end{proposition}

To continue, denote by $A^2_D$ the Hilbert space of holomorphic functions on $D$ whose moduli are square-integrable with respect to the area measure. We equip $A^2_D$ with the norm induced by the inner product
\begin{equation}
\label{innerproduct}
\langle f,g\rangle := \int_Df(z)\overline{g(z)}\,dA.
\end{equation}
Denote by $K_D(z,w)$, $z,w\in D$, the reproducing kernel\footnote{$K_D(z,w)=\frac1\pi\frac{\psi^\prime(z)\overline{\psi^\prime(w)}}{\left(1-\psi(z)\overline{\psi(w)}\right)^2}$, where $\psi$ is any conformal map from $D$ onto $\D$ \cite[\S 1.5]{Gaier}.} for $A^2_D$.  That is,
\begin{equation}
\label{reprkernel}
f(z) = \int_Df(w)K_D(z,w)\,dA
\end{equation}
for any $f\in A^2_D$. It is known \cite[Theorem~I.5.2]{Gaier} that $K_{N,\infty}$ is the reproducing kernel for the set of polynomials of degree at most $N-1$ in the sense of \eqref{reprkernel} and that
\begin{equation}
\label{kernel1}
\big|K_D(z,w) - K_{N,\infty}(z,w) \big| \to 0 \quad \mbox{as} \quad N\to\infty
\end{equation}
locally uniformly for $z,w\in D$.

\begin{theorem}
\label{thm:kernels}
Under the conditions of Theorem~\ref{thm:polynomials}, \eqref{kernel1} holds with $K_{N,\infty}$ replaced by $K_{N,s}$ as $N,s\to\infty$, $N\leq\lfloor s-1\rfloor$.
\end{theorem}

\begin{remark}
The original proof of Theorem~\ref{thm:kernels} as devised by the authors used the full strength of Theorem~\ref{thm:polynomials} and therefore excluded the case $p=0$ and $\ell=1$. The authors are grateful to the anonymous referee who pointed out that only the asymptotics of Carleman polynomials $\pi_{n,\infty}$ on $\overline O$ is needed to prove the theorem and thus allowing all $p\geq0$ and $\ell\in[0,1]$.
\end{remark}

To describe the asymptotic behavior of kernels \eqref{kernel} and \eqref{tildekernel} near the
boundary $\partial K$, it is convenient to introduce the following
notation. Set
\[
H_0(\tau) := 2\frac{e^\tau(\tau-1)+1}{\tau^2} \quad \mbox{and} \quad H_1(\tau) := 6\frac{e^\tau(\tau-2)+\tau+2}{\tau^3},
 \] 
and define $H_{\ell}$ to be the convex combination,
\begin{equation}
\label{H_ell}
H_\ell(\tau) := \frac{3-3\ell}{3-2\ell}H_0(\tau) + \frac{\ell}{3-2\ell} H_1(\tau), \quad \ell\in(0,1).
\end{equation}
Note that the value at the origin for each of these functions is determined by taking a limit; that is, $H_\ell(0)=1$ for all $\ell\in[0,1]$. The following theorem is an analog of \cite[Theorem~2.1]{Lub10}. 

\begin{theorem}
\label{thm:scaling} 
Let $N\leq\lfloor s-1\rfloor$, $z\in\partial K$, $\tau(a,z) := a \Phi'(z) \overline{\Phi(z)}, {\displaystyle \ell := \lim_{N,s \rightarrow  \infty} N s^{-1} \in [0,1]}$ and for
$\ell > 0$, set
\begin{equation}
\label{omegaaz}
\omega(a,z) :=  \left\{
\begin{array}{ll}
\exp\big\{ - \re\left(\tau(a,z)\right)/\ell \big\}, & \re\left(\tau(a,z)\right) > 0, \smallskip \\
1, & \mbox{otherwise}.
\end{array}
\right.
\end{equation}
Under the conditions of Theorem~\ref{thm:polynomials}, assuming $p>0$ when $\ell=1$, it holds that 
\begin{equation}
\label{kernel4}
\lim_{N,s\to\infty}\frac{K_{N,s}(z+\frac aN,z+\frac bN)}{K_{N,s}(z,z)} = H_\ell\left(\tau(a, z) + \overline{\tau(b,z)}\right).
\end{equation}
Moreover, if $\ell>0$, then
\begin{equation}
\label{kernel5}
\lim_{N,s\to\infty}\frac{\widetilde K_{N,s}(z+\frac aN,z+\frac bN)}{\widetilde K_{N,s}(z,z)} = \omega(a, z) \omega(b,z) H_\ell\left(\tau(a, z) + \overline{\tau(b,z)}\right),
\end{equation}
and if $\ell=0$, then
\begin{equation}
\label{eq:9}
\lim_{N,s\to\infty}\frac{\widetilde K_{N,s}(z+\frac aN,z+\frac
  bN)}{\widetilde K_{N,s}(z,z)} = \left\{ 
 \begin{array}{ll}
H_0\left( \tau(a,z) + \overline{\tau(b,z)} \right) &
\quad \re(\tau(a,z)), \re(\tau(b,z)) < 0; \\
0 & \quad \mbox{otherwise.}
  \end{array}
\right.
\end{equation}
The convergence in
(\ref{kernel4})--(\ref{eq:9}) is uniform for $a,b$ in compact subsets
of $\C$ and $z\in\partial K$.   
\end{theorem}

\begin{remark}
The argument of $\tau(a,z)$ is equal to the angle between $a$ and $\map(z)/\map^\prime(z)$, the outward normal to $T$ at $z$.
\end{remark}

\begin{remark}
As is clear from \eqref{tildekernel}, the function  $\omega(a,z)$ is designed to describe the limit of $P_K^{-s}(z+a/N)$  as $N,s\to\infty$. This limit depends on whether or not the points  $z+a/N$ belong to $O$ for $N$ large enough. The case  $\re(\tau(a,z))=0$ corresponds to the situation when the sequence  $\{z+a/N\}$ approaches $z\in\partial K$ tangentially to the  boundary. This does not cause a problem in \eqref{omegaaz} as this  function is continuous with respect to $a$. However, when $\ell=0$  formula \eqref{omegaaz} cannot be used as the limit is described by  a discontinuous function of $a$ and the convexity of the boundary  $\partial K$ at $z$ starts to play a role.
\end{remark}

\begin{remark}
Observe that by putting $s=\infty$ (that is, $\ell=0$), formulae \eqref{kernel2}--\eqref{kernel4} specialize to the asymptotic formulae obtained in \cite{Lub10} for Carleman polynomials. Notice also that when $s=N+1$ ($\ell=1$), the first summands in \eqref{kernel2}--\eqref{kernel4} disappear and only the second ones remain. For general $\ell$, formulae \eqref{kernel2}--\eqref{kernel4} turn out to be convex combination of these two extreme cases. 
\end{remark}

\section{Three Models of Potential Theoretic Ensembles}
\label{sec:two-other-models}

Before proceeding to the proofs of our main results, we will present
three models, a matrix model, an electrostatic model and a polynomial model,
whose joint density of eigenvalues, particles and roots coincide with
the potential theoretic ensembles we are considering.  

\subsection{Entropic Normal Matrix Ensembles}

The {\em entropy} of a self-map $T$ on a metric space $X$ is a
measure of how the distance between nearby points is stretched under
iteration of $T$.  In the case where $\mathbf Z$ is an $N \times N$
complex matrix acting on $\C^N$, the entropy of $\mathbf Z$ is given
by
\[
h(\mathbf Z) = \sum_{n=1}^N \log \max\{1, |\lambda_n| \},
\]
where $\lambda_1, \lambda_2, \ldots, \lambda_N$ are the eigenvalues of
$\mathbf Z$ \cite{springerlink:10.1007/BF01040581}.  We may use this
to create a probability measure on 
normal $N \times N$ complex matrices, which we will denote by
$\mathcal N_N(\C)$.  

There exists a canonical measure on $\mathcal N_N(\C)$ induced by the
standard metric on $\C^{N \times N}$ and we may define a probability
density with respect to this measure by writing
\[
P_N(\mathbf Z) = \frac{1}{Z_N} e^{-2s h(\mathbf Z)},
\]
where $Z_N$ is a normalization constant and $s > N$ is a real number
necessary so that the probability measure is actually finite.  

This probability measure on normal matrices induces a symmetric
probability measure on $\C^N$ as identified with vectors of
eigenvalues.  This measure is absolutely continuous with respect to
Lebesgue measure and its density is given as in (\ref{eq:1}) with
$w(\lambda) = \max\{1, |\lambda|\}^{-2s}$ \cite{MR1643533,
  PhysRevE.55.205}.  Normal matrix ensembles, and in particular the statistics
of their eigenvalues, were first considered in
\cite{LingLieChau|YueYu1992452} and \cite{MR1643533}.  

The function $\lambda \mapsto \log \max\{1, | \lambda | \}$ is the
logarithmic (equilibrium) potential of the closed unit disk, and the
weight for the entropic ensemble is formed from this in the obvious
manner.  We therefore see that the eigenvalue statistics of the entropic
normal matrix ensemble coincides with the potential theoretic ensemble
with $K = \overline{\D}$.

\subsection{Two-Dimensional Electrostatics}

In two-dimensional electrostatics, charged particles are identified
with points in the extended complex plane.  The potential energy of a system of
two like charged particles located at $z, w \in \C$ is proportional to
$-\log| z - w |$.  More generally, if $z_1, z_2, \ldots, z_N$ are the
locations of $N$ identically charged particles, then $\mathbf z$
determines the {\em state} of the system and the potential energy of
this state is given by
\[
 -\sum_{m < n} \log| z_n - z_m |.
\]
The energy is minimized when the particles are all at $\infty$.  In
order for the system to be found in a state where the particles are at
finite positions, there needs to be a potential (or other
obstructions) which repels the particles from $\infty$.  We represent
this field by $V$ so that the interaction energy between a particle
located at $z$ and the field is given by $V(z)$.  The total potential
energy of the system comprised of the $N$ particles in the field is
given by
\[
E(\mathbf z) = \sum_{n=1}^N V(z_n) - \sum_{m < n} \log| z_n - z_m |.
\]

The system is assumed to be in contact with a heat reservoir so that the
energy of the system is variable, but the temperature is fixed.  In
this setting, $\beta$ denotes the reciprocal of the temperature, and
the Boltzmann factor for the state $\mathbf z$ is given by
\[
e^{-\beta E(\mathbf z)} = \bigg\{\prod_{n=1}^N e^{-\beta V(z_n)}
\bigg\} \prod_{m < n} |z_n - z_m|^{\beta}.
\]
This quantity gives the relative density of states, so that the
probability (density) of finding the system in state $\mathbf z$ is
given by
\[
\frac{1}{Z_N} e^{-\beta E(\mathbf z)} \qquad \mbox{where} \qquad Z_N =
\int_{\C^N} e^{-\beta E(\mathbf z)} dA^N(\mathbf z).
\]
Comparing with (\ref{eq:1}) we see that, when $\beta = 2$ the
density of states is identical with the density of eigenvalues of the
normal matrix ensemble with weight $w(z) = e^{-2 V(z)}$.

In this model, a compact set $K$ is identified with a conducting
region.  A charge supported on $K$ will distribute itself to minimize
its potential energy, and this distribution, suitably normalized,
leads to the equilibrium measure on $K$.  In this way, we can think of
the function $-s \log P_K(z)$ as the potential energy felt by an
oppositely charged particle at $z$ when placed in the field given by
the minimal energy configuration formed by placing a total charge of
$s$ on $K$.  In this situation where our system consists of $N$
charged particles, the condition that $s > N$ is required to make
$\infty$ repulsive (or rather to make $K$ sufficiently attractive so
that the particles do not flee to $\infty$).  It follows that the
statistics of particles in this model agree with those of the
potential theoretic ensemble for $K$.  

\subsection{Roots of Random Polynomials}

The {\em Mahler measure} of a polynomial $f(x) \in \C[x]$ is given by
\[
M(f) = \exp \left\{ \int_0^1 \log| f(e^{2 \pi i \theta}) | \, d\theta \right\},
\]
is an example of {\em height function}; that is a function which
measures the complexity of arithmetic objects, in this case
polynomials\footnote{Traditionally Mahler measure is used as a height
  of polynomials in $\mathbb{Q}[x]$, or more generally
  $\overline{\mathbb{Q}}[x]$. 
  However, there is no obstruction in defining it for polynomials in
  $\C[x]$}.  One type of problem of interest 
to number theorists is to provide asymptotic estimates for the number
of arithmetic objects whose height is bounded by $C$ as $C \rightarrow
\infty$.  For instance, for the Mahler measure, such estimates for the
number of integer polynomials of fixed degree and Mahler measure
bounded by $C$ as $C \rightarrow \infty$ was given by Chern and Vaaler
in \cite{MR1868596}.  They also gave a similar estimate for the number of
polynomials with Gaussian integer ($\Z[i]$) coefficients.  

In the latter case, the main term in their estimate came from the
calculation of the Lebesgue measure of the set of polynomials of
degree $N$ with complex coefficients whose Mahler measure is at most
1.  A key aspect of their proof is to show that this volume is equal
to
\[
\frac{\pi}{N+1} \int_{\C^N} \left\{  M\bigg(x^N + \sum_{n=1}^N a_n z^{N-n}
\bigg) \right\}^{-2N-2} \, dA^N(\mathbf a), 
\]
That is the volume is proportional to an integral of a (negative)
power of the Mahler measure of monic polynomials with respect to
Lebesgue measure on the non-leading coefficients of such polynomials.
Moreover, after the change of variables from coefficients to roots of
polynomials, this volume reduces to
\begin{equation}
\label{eq:7}
\frac{\pi}{N+1} \int_{\C^N} \left\{ \prod_{n=1}^N 
\exp \left\{ \int_0^1 \log| \alpha_n - e^{2 \pi i \theta} | \, d\theta \right\}^{-2s}
 \right\} \prod_{m < n} | \alpha_n - \alpha_m |^2   \,
 dA^N(\boldsymbol \alpha); \qquad s = N + 1.
\end{equation}
That is, this volume, up to the factor of $\pi/(N+1)$ is equal to the
normalization constant $Z_N$ for the potential theoretic ensemble for
the unit circle for the value $s = N+1$.  In fact, Chern and Vaaler
were able to show that this normalization constant as a function of
$s$ is a rational function in $s$ with poles at positive integers
$\leq N$.  This striking result can be seen as a consequence of
determinantal nature of the correlation functions.

The derivation of (\ref{eq:7}) shows that the roots of a polynomial
chosen randomly from the volume of complex polynomials of degree $N$
and Mahler measure at most 1 obey the same statistics as those of the
potential theoretic ensemble for the disk.  This gives a polynomial
model for these statistics.  

The computation of the normalization constant of this polynomial model
for potentials for certain other compact regions (in particular the
ellipses considered in Proposition~\ref{prop:ellipse}) is given in
\cite{MR2407817}, while a more general treatment for more general
potentials is given in \cite{MR2707617}.  The special case where the
family of ellipses degenerates to the interval $[-2,2]$ on the real
axis, and its application to the estimation of counting reciprocal
polynomials with bounded Mahler measure is given in \cite{MR2145532}
and \cite{MR2407817}.

\section{Proofs}

To prove Theorem~\ref{thm:polynomials},  we use the method of normal moments in which we rely on the results in \cite[Ch. I]{Suetin}. We must therefore discuss Faber polynomials before proceeding to the proof of Theorem~\ref{thm:polynomials}. We start by stating several auxiliary facts that will be useful later.

\subsection{Auxiliary Facts}
 If $g$ is holomorphic in $\Om$ and vanishes at $\infty$, then for each $r\in[1,\infty)$ the restriction $\overline{g(r\tau)}$, $\tau\in\T:=\{|w|=1\}$,  can be interpreted as the trace on $\T$ of $\overline{g(r/\overline{w})}$, $w\in\D$, which is holomorphic and vanishes at the origin. The latter implies that $\int_\T \tau^k\overline{g(r\tau)}|d\tau|=0$ for all integers $k\geq0$. If, in addition, $h$ is a positive function on $[1,\infty)$ and $\left|w^j\overline{g(w)}h(|w|)\right|$ is integrable with respect to $dA$ for some $j\geq0$, then the Fubini-Tonelli theorem yields that
\begin{equation}
\label{ZeroArea1}
\int_\Om w^j\overline{g(w)}h(|w|)dA = \int_1^\infty \left[\int_\T \tau^j\overline{g(r\tau)}|d\tau|\right]r^{j+1}h(r)dr = 0.
\end{equation}
Furthermore, since $dA(w)=|\map^\prime(z)|^2dA(z)$, where $w=\map(z)$, it holds that
\begin{equation}
\label{ZeroArea2}
\int_O\map^j(z)\map^\prime(z)\overline{G(z)}h(|\map(z)|)dA = 0
\end{equation}
for any $G$ holomorphic in $O$ and vanishing at $\infty$ by \eqref{ZeroArea1} applied with $g=(G\circ\phi)\phi^\prime$, where $\phi$ is the inverse of $\map$ (granted $\left|\map^j\map^\prime Gh(|\map|)\right|$ is integrable with respect to $dA$).

In another connection, the Cauchy-Green identity for the domain $D$ \cite[Thm. 1.2.1]{Hormander} says that
\begin{equation}
\label{CG1} 
\int_D g(z)\overline{h^\prime(z)}\,dA = \frac{1}{2i}\oint_T
g(z)\overline{h(z)}\,dz
\end{equation} 
whenever $g$ and $h^\prime$ are holomorphic functions in $D$ that continuously extend to $T$, where $\oint$ always means integration in the counter-clockwise direction unless specified otherwise. Now, assume that $g$ and $h$ are holomorphic functions in $O$ such that $g$ has at least a double zero at infinity, and $g, h$ and $h^\prime$ continuously extend to $T$. Then by using the transformation  $z\mapsto1/z$ and \eqref{CG1}, one can show that the Cauchy-Green identity for $O$ assumes the form
\begin{equation}
\label{CG2} 
\int_O g(z)\overline{h^\prime(z)}\,dA = -\frac{1}{2i}\oint_T g(z)\overline{h(z)}\,dz.
\end{equation}

\subsection{Faber Polynomials}

Denote by $F_n$ the $n$-th Faber polynomial for $D$ associated with $\map^\prime$. That is,
\[ 
F_n(z) = \oint_T\frac{\map^n(t)\map^\prime(t)}{t-z}\frac{dt}{2\pi i}, \quad z\in D.
\] 
In other words, $F_n$ is the polynomial part of $\map^n\map^\prime$. Then it follows from Plemelj-Sokhotski formulae
\cite{Gakhov} that
\begin{equation}
\label{PlemeljSokhotski} 
F_n = \map^n\map^\prime + E_n \quad \mbox{in} \quad O,
\end{equation} 
where $E_n$ is a holomorphic function vanishing at infinity with integral representation
\begin{equation}
\label{En}
 E_n(z) := \oint_T\frac{\map^n(t)\map^\prime(t)}{t-z}\frac{dt}{2\pi
 i}, \quad z\in O.
\end{equation}

We would like to point out that
\begin{equation}
\label{EnatInfty}
E_n(z) = \mathcal{O}\left(\frac1{z^2}\right) \quad \mbox{as} \quad z\to\infty
\end{equation}
for all integers $n\geq0$. Indeed, consider $\widetilde F_{n+1}$, the
$(n+1)$th Faber polynomial associated with $1$. In this case
\eqref{PlemeljSokhotski} gets replaced by $\widetilde
F_{n+1}=\map^{n+1}+\widetilde E_{n+1}$, where $\widetilde E_{n+1}$ has
an integral representation similar to \eqref{En}. By differentiating
both sides of the last equality, we get that $(n+1)F_n=\widetilde
F_{n+1}^\prime$ and $(n+1)E_n=\widetilde E_{n+1}^\prime$. As
$\widetilde E_{n+1}$ is holomorphic and vanishing at infinity, $E_n$
has at least a double zero there. 

We are interested in the asymptotic behavior of
\begin{equation}
\label{mjks}
m_{k,j}^s := \int_\C F_j\overline{F_k}P_K^{-2s}\,dA = \int_D F_j\overline{F_k}\,dA + \int_O F_j\overline{F_k}|\map|^{-2s}\,dA,
\quad j,k\leq \lfloor s-2 \rfloor, 
\end{equation}
where we used \eqref{anotherform} for the second representation. 

It was shown in \cite[Equation (1.32) combined with \eqref{CG2} above]{Suetin} that the first integral on the right-hand side of \eqref{mjks} can be written as 
\begin{equation}
\label{ID}
\frac{\pi}{k+1}\left(\delta_{jk}-\frac{k+1}{\pi}\int_OE_j\overline{E_k}\,dA\right) =:\frac{\pi}{k+1}\left(\delta_{kj}+I_D\right),
\end{equation}
where the integral over $O$ is well defined (finite) by \eqref{EnatInfty}. Moreover,  it was also obtained there, see \cite[Equation (1.45) and Lemma~1.5]{Suetin}, that
\begin{equation}
\label{IDbound}
|I_D| \leq \frac{\const}{(j+1)^{p+\alpha}(k+1)^{p+\alpha}} \quad \mbox{or} \quad |I_D| \leq \const \rho^{j+k},
\end{equation}
where both constants are independent of $j$ and $k$, but depend on $T$
and $\rho$ (in the analytic case). Hereafter, by stating a double
estimate of the form \eqref{IDbound}, we always assume that the first
bound is given for $T$ of class $C^{p+1,\alpha}$ and the second one
for $T$ analytic with $\rho(T)<\rho<1$.  

On the other hand, the second integral on the right-hand side of \eqref{mjks} can be written with the help of \eqref{PlemeljSokhotski} as
\begin{equation}
\label{fourint}
\int_O\map^j\overline{\map^k}|\map^\prime|^2|\map|^{-2s}dA + \int_O\map^j\map^\prime\overline{E_k}|\map|^{-2s}dA  + \int_O E_j\overline{\map^k\map^\prime}|\map|^{-2s}dA + \int_O E_j\overline{E_k}|\map|^{-2s}dA.
\end{equation}
It can be immediately computed by conformality of $\map$ that the first integral in \eqref{fourint} is equal to
\[
\int_O\map^j\overline{\map^k}|\map^\prime|^2|\map|^{-2s}dA = \int_\Om w^j\overline{w^k}|w|^{-2s}dA = \frac{\pi}{s-(k+1)}\delta_{kj}.
\]
The second integral in \eqref{fourint} is zero by \eqref{ZeroArea2} applied with $G=E_k$ and $h(r)=r^{-2s}$, $r\in[1,\infty)$. Conjugating the third integral in \eqref{fourint}, analogous reasoning shows it is zero as well. Thus, similar to \eqref{ID}, the second integral on the right-hand side of \eqref{mjks} can be written as
\begin{equation}
\label{IO}
\frac{\pi}{s-(k+1)}\left(\delta_{jk}+\frac{s-(k+1)}{\pi}\int_OE_j\overline{E_k}|\map|^{-2s}dA\right) =: \frac{\pi}{s-(k+1)}\left(\delta_{kj}+I_O\right).
\end{equation}

We claim that $|I_O|$ satisfies \eqref{IDbound} as well, namely,
\begin{equation}
\label{IObound}
|I_O| \leq \frac{\const}{(j+1)^{p+\alpha}(k+1)^{p+\alpha}} \quad \mbox{or} \quad |I_O| \leq \const \rho^{j+k}.
\end{equation}
Indeed, to prove \eqref{IDbound}, it was shown in \cite[(1.44) and the following paragraph, (1.45)]{Suetin} that
\begin{equation}
\label{boundEk}
\int_\T\big|(E_n\circ\phi)(\tau)\phi^\prime(\tau)\big|^2|d\tau| \leq  \frac{\const}{(n+1)^{2p+2\alpha}} \quad \mbox{or} \quad  \leq \const \rho^{2n}
\end{equation} 
for all integers $n\geq0$, where $\phi$ is the inverse of $\map$. Moreover, the monotonicity of $L^2$-norms implies that if the integrand in \eqref{boundEk} is evaluated at $r\tau$ rather than at $\tau$ for any fixed $r\in(1,\infty)$, the estimate remains valid. Then on account of
\[
\left|\int_OE_j\overline{E_k}|\map|^{-2s}dA\right| \leq \int_1^\infty\left[\int_\T\left|\big((E_jE_k)\circ\phi\big)(r\tau)\right|\big|\phi^\prime(r\tau)\big|^2|d\tau|\right]r^{1-2s}dr
\]
and the Cauchy-Schwarz inequality, $|I_O|$ is bounded by
\[
\frac{s-(k+1)}{\pi}\left(\int_1^\infty r^{1-2s}dr\right)\left(\int_\T\big|(E_k\circ\phi)(\tau)\phi^\prime(\tau)\big|^2|d\tau|\right)^{1/2}\left(\int_\T\big|(E_j\circ\phi)(\tau)\phi^\prime(\tau)\big|^2|d\tau|\right)^{1/2}.
\]
Clearly, \eqref{IObound} follows now from \eqref{boundEk}.

Finally, gathering together \eqref{ID} and \eqref{IO}, we get that
\begin{equation}
\label{mjks2}
m_{k,j}^s = \frac{s\pi}{(k+1)(s-(k+1))}\left(\delta_{kj}+\epsilon_{k,j}^s\right),
\end{equation}
where
\begin{equation}
\label{epsilonjks}
\epsilon_{k,j}^s := \frac{s-(k+1)}{s}I_D+\frac{k+1}{s} I_O = -\frac{k+1}{\pi}\left(1-\frac{k+1}{s}\right)\int_OE_j\overline{E_k}\left(1-|\map|^{-2s}\right)\,dA
\end{equation}
and
\begin{equation}
\label{estepsilon}
|\epsilon_{k,j}^s| \leq \frac{\const}{(j+1)^{p+\alpha}(k+1)^{p+\alpha}} \quad \mbox{or} \quad |\epsilon_{k,j}^s| \leq \const \rho^{j+k}
\end{equation}
by \eqref{IDbound} and \eqref{IObound}.

\subsection{The Von Koch-Riesz Algebra}

Denote by $\mathcal{D}$ the algebra of all operators defined on $\ell_2(\N)$ by matrices $\mathbf{A}=[a_{k,j}]_{j,k=0}^{\infty}$ with respect to the standard basis for which
\[
\|\mathbf{A}\|_{\mathcal D} := \max\left\{\sum_{k=0}^\infty|a_{k,k}|,\left(\sum_{j=0}^\infty\sum_{k=0}^\infty|a_{k,j}|^2\right)^{1/2}\right\} <\infty.
\]
It is known \cite[Theorem~II.2.1]{GohbergGoldbergKrupnik} that if $\{\mathbf A_n\}$ is a sequence in $\mathcal{D}$ converging to $\mathbf A \in \mathcal D$ (with respect to $\| \cdot \|_{\mathcal D}$), then the determinant of $\mathbf I + \mathbf A_n$ ($\mathbf I$ being the identity operator) converges to the determinant of $\mathbf I +\mathbf A$.  

Let $\{s_n\}$ be an increasing sequence of positive reals such that $s_n\to\infty$ as $n\to\infty$. Set, for convenience, $\epsilon_{k,j}^{s_n}:=0$ when either $j$ or $k$ is greater than $\lfloor s_n-2\rfloor$ and define $\mathbf{E}_{s_n}:=[\epsilon_{k,j}^{s_n}]_{j,k=0}^\infty$. We also set $\mathbf{E}_\infty:=[\epsilon_{k,j}^\infty]_{j,k=0}^\infty$, where we put
\begin{equation}
\label{epsilonjkinf}
\epsilon_{k,j}^\infty := -\frac{k+1}{\pi}\int_OE_j\overline{E_k}\,dA.
\end{equation}
Observe that the estimate in \eqref{estepsilon} is also valid for $s=\infty$. Using this bound, it is simple to verify that
\[
\|\mathbf{E}_{s_n}\|_\mathcal{D} \leq \const \sum_{k=1}^\infty\frac{1}{k^{2(p+\alpha)}} \quad \mbox{or} \quad \|\mathbf{E}_{s_n}\|_\mathcal{D} \leq \frac{\const}{1-\rho^2}
\]
for each $n$ including the case $n=\infty$, where the constant $\sum_{k=1}^\infty k^{-2(p+\alpha)}$ is finite as $p+\alpha>1/2$. Thus, all the operators $\mathbf{E}_{s_n}$ belong to the Von Koch-Riesz algebra $\mathcal{D}$. Moreover, it holds that
\begin{equation}
\label{normto0}
\|\mathbf{E}_{s_n}-\mathbf{E}_\infty\|_\mathcal{D} \to 0 \quad \mbox{as} \quad n\to\infty.
\end{equation}
Indeed, let $\{k_n\}$ be a non-decreasing sequence of integers such that $k_n\to\infty$ and $k_n/n\to 0$ as $n\to\infty$. Then
\begin{equation}
\label{sumkN1}
\sum_{k=k_n}^\infty\left|\epsilon_{k,k}^{s_n}-\epsilon_{k,k}^\infty\right| \leq \sum_{k=k_n}^\infty\left(|\epsilon_{k,k}^{s_n}|+|\epsilon_{k,k}^\infty|\right) \leq \frac{\const}{(k_n+1)^{p+\beta}} \quad \mbox{or} \quad \sum_{k=k_n}^\infty\left|\epsilon_{k,k}^{s_n}-\epsilon_{k,k}^\infty\right| \leq \const\rho^{2k_n}
\end{equation}
by \eqref{estepsilon}. Furthermore, we can readily deduce from \eqref{epsilonjks} and \eqref{epsilonjkinf} using the notation of \eqref{ID} and \eqref{IO} that
\[
\left|\epsilon_{k,k}^{s_n}-\epsilon_{k,k}^\infty\right| = \frac{k+1}{n}\left|I_O-I_D\right| \leq \frac{k_n}{n}\frac{\const}{(k+1)^{2(p+\alpha)}} \quad \mbox{or} \quad \left|\epsilon_{k,k}^{s_n}-\epsilon_{k,k}^\infty\right|\leq \frac{k_n}{n}\const\rho^{2k}
\]
by \eqref{IDbound} and \eqref{IObound} for all $k\in\{0,\ldots,k_n-1\}$. Therefore, it holds that
\begin{equation}
\label{sumkN2}
\sum_{k=0}^{k_n-1}\left|\epsilon_{k,k}^{s_n}-\epsilon_{k,k}^\infty\right| \leq \const\frac{k_n}{n}.
\end{equation}
Combining \eqref{sumkN1} and \eqref{sumkN2}, we deduce that
\begin{equation}
\label{fullsum1}
\sum_{k=0}^\infty\left|\epsilon_{k,k}^{s_n}-\epsilon_{k,k}^\infty\right| \to 0 \quad \mbox{as} \quad n\to\infty
\end{equation}
by the choice of the sequence $\{k_n\}$. Analogously, one can show that
\begin{equation}
\label{fullsum2}
\sum_{j=0}^\infty\sum_{k=0}^\infty\left|\epsilon_{k,j}^{s_n}-\epsilon_{k,j}^\infty\right|^2 \to 0 \quad \mbox{as} \quad n\to\infty,
\end{equation}
which finishes the proof of \eqref{normto0}.

Naturally \cite[Section~I.1]{GohbergGoldbergKrupnik}, it holds that $\det(\mathbf{I}+\mathbf{E}_{s_n}) = \det\left[\delta_{kj}+\epsilon_{k,j}^{s_n}\right]_{j,k=0}^{\lfloor s_n-2\rfloor}$ and therefore we deduce from the remark made at the beginning of this section that
\begin{equation}
\label{limitdet}
\det(\mathbf{I}+\mathbf{E}_{s_n}) \to \det(\mathbf{I}+\mathbf{E}_\infty)>0 \quad \mbox{as} \quad n\to\infty,
\end{equation}
where the last inequality was shown in \cite[Section~I.4]{Suetin}.

\subsection{Proof of Theorem~\ref{thm:polynomials}}

Since $\{F_n\}$ is a complete system of polynomials, each $\pi_{n,s}$
can be expressed as a linear combination of $F_0,\ldots,F_n$ with the
coefficients determined via the orthogonality relations
\eqref{orthogonality}. In fact, it holds that 
\begin{equation}
\label{pins}
\pi_{n,s}(z) = \frac{1}{\sqrt{D_{n-1,s}D_{n,s}}}\left[
\begin{array}{cccc} 
m_{0,0}^s & m_{0,1}^s & \ldots & m_{0,n}^s \\
m_{1,0}^s & m_{1,1}^s & \ldots & m_{1,n}^s \\ 
\vdots & \vdots & \ddots & \vdots \\ 
m_{n-1,0}^s & m_{n-1,1}^s & \ldots & m_{n-1,n}^s \\ 
F_0(z) & F_1(z) & \ldots & F_n(z)
\end{array} 
\right],
\end{equation}
where the moments $m_{k,j}^s$ are defined in \eqref{mjks} and $D_{n,s}:=\det[m_{k,j}^s]_{j,k=0}^n$. 

Set $\Delta_{n,s}:=\det[\delta_{kj}+\epsilon_{k,j}^s]_{j,k=0}^n$ and observe that
\begin{equation}
\label{Ddelta}
D_{n,s} = \Delta_{n,s}\prod_{k=0}^n\frac{s\pi}{(k+1)(s-(k+1))} \quad \mbox{and} \quad D_{n,s}(j) = \Delta_{n,s}(j)\prod_{k=0,\,k\neq j}^n\frac{s\pi}{(k+1)(s-(k+1))}
\end{equation}
by \eqref{mjks2}, where the determinants $D_{n,s}(j)$ and $\Delta_{n,s}(j)$ are obtained from the same matrices as $D_{n,s}$ and $\Delta_{n,s}$ only with the last row and the $(j+1)$-st column removed. Given \eqref{estepsilon}, it is a straightforward algebraic computation using Hadamard's inequality, see \cite[Lemma~1.7]{Suetin}, to derive that
\begin{equation}
\label{Deltanskest}
\Delta_{n,s}(j) \leq \frac{\const}{(j+1)^{p+\beta}(n+1)^{p+\beta}} \quad \mbox{or} \quad \Delta_{n,s}(j) \leq \const \rho^{j+n}
\end{equation}
for any $j\in\{0,1,\ldots,n-1\}$. 

On the other hand, the family $\{\Delta_{n,s}\}$ is bounded away from zero. Indeed, as mentioned just before \eqref{anotherform}, $\map^\prime(\infty)=\gamma_K^{-1}$ and therefore the leading coefficient of $F_n$ is equal to $\gamma_K^{-n-1}$. Hence, we get from \eqref{pins} and \eqref{Ddelta} that
\begin{equation}
\label{kappadelta}
\varkappa_{n,s}\gamma_K^{n+1} = \sqrt{\frac{D_{n-1,s}}{D_{n,s}}} = \sqrt{\frac{n+1}{\pi}\left(1-\frac{n+1}{s}\right)\frac{\Delta_{n-1,s}}{\Delta_{n,s}}}.
\end{equation}
Recall that any monic orthogonal polynomial has the smallest
$L^2$-norm with respect to the weight of orthogonality among all monic
polynomials of the same degree. In particular, 
\[
\frac{1}{\varkappa_{n,s}^2} = \int_\C|\pi_{n,s}/\varkappa_{n,s}|^2|P_K|^{-2s}dA \leq \int_\C|\gamma_K^{n+1}F_n|^2|P_K|^{-2s}dA = \gamma_K^{2n+2}m_{n,n}^s.
\]
Therefore, it follows from \eqref{kappadelta}, \eqref{mjks2}, and \eqref{epsilonjks} that
\begin{eqnarray}
\Delta_{n,s} &=& \frac{n+1}{\pi}\left(1-\frac{n+1}{s}\right)\frac{\Delta_{n-1,s}}{\gamma_K^{2n+2}\varkappa_{n,s}^2} \leq \frac{n+1}{\pi}\left(1-\frac{n+1}{s}\right) m_{n,n}^s\Delta_{n-1,s} \nonumber \\
&=& (1+\epsilon_{n,n}^s)\Delta_{n-1,s}< \Delta_{n-1,s}. \nonumber
\end{eqnarray}
Hence, it holds that
\begin{equation}
\label{infdeltans}
\inf_s\min_{1\leq n\leq \lfloor s-2\rfloor}\Delta_{n,s} = \inf_s\Delta_{\lfloor s-2\rfloor,s}>0
\end{equation}
by \eqref{limitdet} since $\Delta_{\lfloor s-2\rfloor,s}=\det(\mathbf{I}+\mathbf{E}_s)$, which proves the claim.

Thus, expanding the determinant $\Delta_{n,s}$ by the last row, we get that
\begin{equation}
\label{expanddelta}
\Delta_{n,s} = (1+\epsilon_{n,n}^s)\Delta_{n-1,s} + \sum_{j=0}^{n-1}(-1)^{n+j}\epsilon_{n,j}^s\Delta_{n,s}(j).
\end{equation}
Dividing both sides of the equality above by $\Delta_{n-1,s}$ and using \eqref{estepsilon}, \eqref{Deltanskest}, and \eqref{infdeltans} yields
\begin{equation}
\label{ratiodelta}
\frac{\Delta_{n,s}}{\Delta_{n-1,s}} = 1 + \mathcal{O}\left(\frac{1}{n^{2(p+\alpha)}}\right) \quad \mbox{or} \quad \frac{\Delta_{n,s}}{\Delta_{n-1,s}} = 1 + \mathcal{O}\left(\rho^{2n}\right).
\end{equation}
Clearly, we get \eqref{kappa} by taking the reciprocal of \eqref{ratiodelta} and substituting it into \eqref{kappadelta}.

Now, expanding the determinant in \eqref{pins} by the last row as in \eqref{expanddelta} yields
\[
\pi_{n,s} = \sqrt{\frac{n+1}{\pi}\left(1-\frac{n+1}{s}\right)}\sqrt\frac{\Delta_{n-1,s}}{\Delta_{n,s}}\left(F_n+\sum_{j=0}^{n-1}(-1)^{n+j}\frac{(j+1)(s-j-1)}{(n+1)(s-n-1)}\frac{\Delta_{n,s}(j)}{\Delta_{n-1,s}}F_j\right).
\]
Hence, by factoring out $\map^n\map^\prime$ and using \eqref{PlemeljSokhotski}, the error term in \eqref{asymp} can be written as
\[
\sqrt\frac{\Delta_{n-1,s}}{\Delta_{n,s}}-1+\sqrt\frac{\Delta_{n-1,s}}{\Delta_{n,s}}\left[\frac{E_n}{\map^n\map^\prime}+\sum_{j=0}^{n-1}(-1)^{n+j}\frac{(j+1)(1-\frac{j+1}{s})}{(n+1)(1-\frac{n+1}{s})}\frac{\Delta_{n,s}(j)}{\Delta_{n-1,s}}\left(1+\frac{E_j}{\map^j\map^\prime}\right)\frac{1}{\map^{n-j}}\right].
\]
Since $|\map| > 1$ and $|\map^\prime|$ is bounded away from zero in $O$, we get from \eqref{boundEk}, \eqref{Deltanskest}, \eqref{infdeltans}, and \eqref{ratiodelta} that the error term in \eqref{asymp} is of order 
\begin{equation}
\label{Xns}
\frac{\log(n+1)}{(n+1)^{p+\alpha}} + \frac{1}{(n+1)^{p+\alpha}}\sum_{j=0}^{n-1}\frac{(j+1)(1-\frac{j+1}{s})}{(n+1)(1-\frac{n+1}{s})}\frac{1}{(j+1)^{p+\alpha}}.
\end{equation}
If $\limsup_{n,s\to\infty} n/s<1$, then the fractions $(1-\frac{j+1}{s})/(1-\frac{n+1}{s})$ are uniformly bounded above and it easily follows from \eqref{Xns} that the error term in \eqref{asymp} is of order
\begin{equation}
\label{Xns1}
\frac{\log(n+1)}{(n+1)^{p+\alpha}} + \frac{1}{(n+1)^{2(p+\alpha)-1}},
\end{equation}
where the first summand is larger for all $p\geq1$ and the second one is larger when $p=0$. Clearly, the estimate for the error term in the case of analytic curve, can be derived in a similar fashion. On the other hand, if $\limsup_{n,s\to\infty} n/s=1$, then we use the estimate
\[
1-\frac{n+1}{s} \geq 1-\frac{n+1}{n+2} = \frac{1}{n+2},
\]
which is valid since $n\leq\lfloor s-2\rfloor$. In this case \eqref{Xns1} gets replaced by
\[
\frac{\log(n+1)}{(n+1)^{p+\alpha}} + \frac{1}{(n+1)^{2(p+\alpha)-2}},
\]
where the first summand is larger for all $p\geq2$ and and the second one is larger when $p=1$ (we exclude $p=0$ as in this case the above bound grows as $n^{2(1-\alpha)}$). Analogous estimate shows that the error term in \eqref{asymp} is of order $\rho^{-n}$ when $T$ is an analytic curve. This finishes the proof of Theorem~\ref{thm:polynomials}. \qed

\subsection{Proof of Proposition~\ref{prop:ellipse}} Let $U_n$ be the monic Chebysh\"ev polynomial of the second kind for the interval $\left[-2\sqrt q,2\sqrt q\right]$. That is,
\[ 
U_n = \map^n\map^\prime\left(1-\frac{q^{n+1}}{\map^{2n+2}}\right), \quad \map(z) = \frac{z+\sqrt{z^2-4q}}{2}, \quad z\in O.
\]
It can be readily checked that the inverse of $\map$ is indeed $\phi(w)=w+q/w$, $\map$ is the conformal map of the complement of $\left[-2\sqrt q,2\sqrt q\right]$ onto $\left\{w:~|w|>\sqrt q\right\}$ with positive derivative at infinity, and the level lines of $\map$ are ellipses with foci $\pm2\sqrt q$.

Let us show that polynomials $U_n$ are orthogonal on $D$ with respect
to area measure. It follows from the Cauchy-Green identity \eqref{CG1}
that 
\begin{eqnarray}
 2i(k+1)\int_DU_n(z)\overline{z^k}dA &=& \oint_T \left(\map^n(z)-q^{n+1}/\map^{n+2}(z)\right)\map^\prime(z)\overline{z^{k+1}}dz \nonumber
\label{intD}
\\ &=& \oint_\T\left(\tau^n-q^{n+1}/\tau^{n+2}\right)(\tau q+1/\tau)^{k+1}d\tau = 0
\end{eqnarray} 
for $k<n$, where we used the identity $\tau=1/\bar \tau$ on $\T$ and the last equality is a consequence of the facts $\oint_\T \tau^n\tau^jd\tau=0$ for all $j\geq -n$ and $\oint_\T \tau^{-n-2}\tau^jd\tau=0$ for all $j\leq n$. 

In another connection, it holds that
\begin{eqnarray}
\int_O U_n(z)\overline{z^k}|\map(z)|^{-2s}dA &=& \int_O\big(\map^n(z)-q^{n+1}/\map^{n+2}(z)\big)\overline{\big(z^k/\map^\prime(z)\big)}|\map(z)|^{-2s}|\map^\prime(z)|^2dA \nonumber \\
&=& \int_\Om\big(w^n-q^{n+1}/w^{n+2}\big)\overline{\phi^k(w)\phi^\prime(w)}|w|^{-2s}dA. \nonumber
\end{eqnarray}
It is easy to check using the expressions $\phi(w)=w+q/w$ and $w=r\tau$, $r\in[1,\infty)$ and $\tau\in\T$, that the chain of equalities above can be continued as
 \begin{equation}
 \label{intO}
\int_1^\infty\left[\int_\T\left((r\tau)^n-\frac{q^{n+1}}{(r\tau)^{n+2}}\right)\left(\frac r\tau+\frac{q\tau}{r}\right)^k\left(1-\frac{q\tau^2}{r^2}\right)|d\tau|\right]\frac{rdr}{r^{2s}} = 0
\end{equation}
for $k<n$, since the Laurent polynomials in $\tau$ integrated over $\T$ does not contain a constant term as the exponents of $\tau$ range from $n+k+2$ down to $n-k$ and then from $k-n$ down to $-k-n$.

Altogether, the polynomials $U_n$ are orthogonal over $\C$ with
respect to the measure $P_K^{-2s}dA$. In fact, it can be easily shown
that they are also the Faber polynomials for this $K$. It remains to
compute the normalizing factor. Evaluating as in 
\eqref{intD} and \eqref{intO}, we get that
\begin{eqnarray}
\int_D U_n(z)\overline{z^n}dA &=& \frac{1}{2i}\frac{1}{n+1}\oint_\T \left(w^n-\frac{q^{n+1}}{w^{n+2}}\right)\left((qw)^{n+1}+\cdots+\frac{1}{w^{n+1}}\right)dw \nonumber \\
&=& \frac{\pi}{n+1}\left(1-q^{2n+2}\right) \nonumber
\end{eqnarray}
and
\begin{eqnarray}
\int_O U_n(z)\overline{z^n}|\map(z)|^{-2s}dA &=& \int_1^\infty\left[\int_\T\left((r\tau)^n-\frac{q^{n+1}}{(r\tau)^{n+2}}\right)\left(-\frac{q^{n+1}\tau^{n+2}}{r^{n+2}}+\cdots+\frac{r^n}{\tau^n}\right)|d\tau|\right]\frac{rdr}{r^{2s}} \nonumber \\
&=& \pi\left(\frac{1}{s-(n+1)}+\frac{q^{2n+2}}{s+n+1}\right). \nonumber
\end{eqnarray}
Thus, we deduce that
\begin{eqnarray}
\varkappa_{n,s}^{-2}&=&\int_\C |U_n(z)|^2P_K^{-2s}(z)dA = \int_D U_n(z)\overline{z^n}dA + \int_O U_n(z)\overline{z^n}|\map(z)|^{-2s}dA \nonumber \\
&=& \frac{\pi s}{(n+1)(s-n-1)}\left(1 - q^{2n+2}\frac{s-n-1}{s+n+1}\right). \nonumber \qed
\end{eqnarray}

\subsection{Proof of Proposition~\ref{prop:kernels}}

Fix $N\leq\lfloor s-1\rfloor$ and let $\zeta$ be a point such that
\begin{equation}
\label{notinK}
\zeta\in\overline O \quad \mbox{and} \quad \dist(\zeta,K)\leq c/N
\end{equation}
for some fixed constant $c$. Further, let $\zeta_0\in T$ be such that
$|\zeta-\zeta_0|=\dist(\zeta,K)$. Since $|\map(\zeta_0)|=1$ and $\map$
is continuously differentiable in $\overline O$ (since $T$ is at least
$C^{1,\alpha}$-smooth) it holds that 
\[
|\map(\zeta)| \leq 1 + |\map(\zeta)-\map(\zeta_0)| \leq 1 + \mathcal{O}(|\zeta-\zeta_0|) = 1 + \mathcal{O}(N^{-1}),
\]
where the estimate $\mathcal{O}(\cdot)$ does not depend on the choice of $\zeta$ satisfying \eqref{notinK}. Hence,
\begin{equation}
\label{maxvarphik}
\max_{k\in\{1,\ldots,N\}}|\map^k(\zeta)| \leq \const
\end{equation}
for some absolute constant. 

Select $z,w$ satisfying \eqref{notinK} and assume that $z\neq w$. Put, for brevity, $u:=\map(z)\overline{\map(w)}$. Then we get from the definition of $K_{N,s}$, \eqref{asymp}, and \eqref{maxvarphik} that
\begin{eqnarray}
K_{N,s}(z,w) &=& \frac{\map^\prime(z)\overline{\map^\prime(w)}}{\pi}\sum_{n=0}^{N-1}\left((n+1)-\frac{(n+1)^2}{s}\right)u^n\big[1+\mathcal{O}\left(\Sigma_n\right)\big]\nonumber \\
&=& \frac{\map^\prime(z)\overline{\map^\prime(w)}}{\pi}\left(\sum_{n=0}^{N-1}(n+1)u^n-\frac1s\sum_{n=0}^{N-1}(n+1)^2u^n\right) \nonumber\\
\label{KNsSum}
&&+ \mathcal{O}\left(\max\left\{1,N^2\Sigma_N\right\}\right).
\end{eqnarray}
Since,
\[
\sum_{n=0}^{N-1}(n+1)u^n =-(N+1)\frac{u^N}{1-u}+\frac{1-u^{N+1}}{(1-u)^2},
\]
and
\[
\sum_{n=0}^{N-1}(n+1)^2u^n = -(N+1)^2\frac{u^N}{1-u} + (N+1)\frac{1-u^{N+1}}{(1-u)^2}- (N+2)\frac{1+u^{N+1}}{(1-u)^2} + 2\frac{1-u^{N+2}}{(1-u)^3},
\]
the validity of \eqref{kernel2} follows. In a similar but simpler fashion, we also get \eqref{kernel3}.

\subsection{Proof of Theorem~\ref{thm:kernels}}

Recall the Christoffel variational principle:
\begin{equation}
\label{Christoffel}
K_{N,s}(z,z) = \max_{\deg(p)<N}\frac{|p(z)|^2}{\int_\C|p|^2P_K^{-2s}dA}, \quad z \in \C,
\end{equation}
and the reproducing property:
\begin{equation}
\label{reprproperty}
p(z) = \int_\C p(w)K_{N,s}(z,w)P_K^{-2s}(w)dA, \quad \deg(p)<N,
\end{equation}
that hold for all $s\in(1,\infty]$ and $N\leq\lfloor s-1\rfloor$. It can be readily deduced from \eqref{Christoffel} that
\begin{equation}
\label{chain}
K_{N,s}(z,z)\leq K_{N,\infty}(z,z)\leq K_D(z,z),
\end{equation}
where the second inequality follows from the fact that
$K_D(z,w)=\sum_{n=0}^\infty
\overline{\pi_{n,\infty}(w)}\pi_{n,\infty}(z)$
\cite[\S~1.5]{Gaier}. Furthermore, \eqref{reprproperty} together with
\eqref{chain} yield  
\begin{equation}
\label{reprDO}
\int_D|K_{N,s}(z,w)|^2dA\leq K_D(z,z) \quad \mbox{and} \quad \int_O|K_{N,s}(z,w)||\map(w)|^{-2s}dA \leq K_D(z,z).
\end{equation}

It follows from \eqref{reprkernel} and \eqref{reprDO} that
\begin{eqnarray}
\int_D|K_D(u,w)-K_{N,s}(u,w)|^2dA &=& K_D(w,w)-2K_{N,s}(w,w)+\int_D|K_{N,s}(u,w)|^2dA\nonumber\\
& \leq & K_D(w,w)-K_{N,s}(w,w).\nonumber
\end{eqnarray}
Therefore, \eqref{reprkernel}, the Cauchy-Schwarz inequality, and the above estimate yield 
\begin{eqnarray}
\left|K_D(z,w)-K_{N,s}(z,w)\right| & \leq & \int_D\left|K_D(u,w)-K_{N,s}(u,w)\right||K_D(z,u)|dA\nonumber \\
&\leq& \big(K_D(z,z)\big)^{1/2}\big(K_D(w,w)-K_{N,s}(w,w)\big)^{1/2}. \nonumber
\end{eqnarray}
That is, we only need to demonstrate the convergence in
\eqref{kernel1} along the diagonal. Moreover, since \eqref{kernel1} is
valid for $K_{N,\infty}$, it suffices to show only that
$K_{N,\infty}(z,z) - K_{N,s}(z,z) \to 0$ locally uniformly in $D$ as
$N,s\to\infty$. To this end, observe that 
\begin{eqnarray}
K_{N,\infty}(z,z) &=& \int_\C K_{N,\infty}(z,u)K_{N,s}(u,z)P_K^{-2s}(u)dA\nonumber\\
&=& K_{N,s}(z,z) + \int_O K_{N,\infty}(z,u)K_{N,s}(u,z)|\map(u)|^{-2s}dA \nonumber
\end{eqnarray}
by \eqref{reprproperty} and the positivity of $K_{N,\infty}(z,z)$. Then we get from the equality above, the Cauchy-Schwarz inequality, and \eqref{reprDO} that
\begin{eqnarray}
K_{N,\infty}(z,z)- K_{N,s}(z,z) &\leq& \left(\int_O\frac{|K_{N,s}(z,u)|^2}{|\map(u)|^{2s}}dA\right)^{1/2}\left(\int_O\frac{|K_{N,\infty}(z,u)|^2}{|\map(u)|^{2s}}dA\right)^{1/2}\nonumber\\
\label{referee}
&\leq& \big(K_D(z,z)\big)^{1/2}\left(\int_O\frac{|K_{N,\infty}(z,u)|^2}{|\map(u)|^{2s}}dA\right)^{1/2}.
\end{eqnarray}

To estimate the integral in \eqref{referee}, observe that
\[
|K_{N,\infty}(z,u)|^2 \leq K_{N,\infty}(z,z)K_{N,\infty}(u,u) \leq K_D(z,z)K_{N,\infty}(u,u)
\]
by the Cauchy-Schwarz inequality and \eqref{chain}. Hence,
\begin{eqnarray}
\int_O|K_{N,\infty}(z,u)|^2|\map(u)|^{-2s}dA &\leq& K_D(z,z)\sum_{n=0}^{N-1}\int_O|\pi_{n,\infty}|^2|\map|^{-2s}dA\nonumber\\
&\leq& \const K_D(z,z)\sum_{n=0}^{N-1}\int_O|\map|^{2n-2s}|\map^\prime|^2dA\nonumber\\
&\leq& \const K_D(z,z)\sum_{n=0}^{N-1}\frac{n+1}{s-(n+1)} \nonumber\\
\label{almostthere}
&=& \const K_D(z,z)\frac{N^2}{s-N},
\end{eqnarray}
where we used \eqref{asymp} for $s=\infty$. Clearly,
\eqref{almostthere}, \eqref{referee}, and the reduction process
carried out above, prove Theorem~\ref{thm:kernels} under the condition
$N^2/s\to0$ as $N,s\to\infty$. The proof for all $N$ readily follows
from the obvious inequlity $K_{N,s}(z,z)\leq K_{M,s}(z,z)$, $N\leq
M$. 
\qed

\subsection{Proof of Theorem~\ref{thm:scaling}}

Since $p>0$ when $\ell=1$, it holds that $\Sigma_N\to0$ as $N\to\infty$. Thus, we deduce from \eqref{kernel3} that
\begin{equation}
\label{Limit1}
\lim_{N,s\to\infty} K_{N,s}(z,z)N^{-2} = \frac{|\map^\prime(z)|^2}{\pi}\left(\frac{1-\ell}{2} +  \frac{\ell}{6}\right) =\frac{|\map^\prime(z)|^2}{\pi}\frac{3-2\ell}{6}
\end{equation} 
uniformly for $z\in T$. Fix $c>0$ and let $a,b$ be such that $|a|,|b|<c$.  It follows from the Cauchy-Schwarz inequality that 
\begin{equation}
\label{bound1}
\left|K_{N,s}\left(z+\frac aN,z+\frac bN\right)\right| \leq \max|K_{N,s}(\zeta,\zeta)|,
\end{equation}
where the maximum is taken over all $\zeta$ satisfying \eqref{notinK}. Furthermore, we get from the Bernstein-Walsh inequality and \eqref{maxvarphik} that
\begin{equation}
\label{bound2}
|K_{N,s}(\zeta,\zeta)| \leq |\map(\zeta)|^{2(N-1)}\max_{w\in T}|K_{N,s}(w,w)| \leq\const\max_{w\in T}|K_{N,s}(w,w)|
\end{equation}
for any $\zeta$ satisfying \eqref{notinK} with some absolute constant. Combining \eqref{bound1} and \eqref{bound2} with \eqref{Limit1}, we see that $\left\{K_{N,s}\left(z+\frac aN,z+\frac bN\right)N^{-2}\right\}$ is a normal family for $|a|,|b|<c$, where the functions in this family are indexed by $z\in T$, $N\in\N$, and $s\in[N+1,\infty)$. Therefore, it suffices to prove \eqref{kernel4} only for those $a,b$ for which $z+\frac aN,z+\frac bN\notin K$ (since $T$ has a tangent at $z$, it holds that either $z+a/N\in K$ for all $N$ large enough or $z+a/N\not\in K$ for all $N$ large enough). 

As $|\map^\prime|$ is bounded above in $\overline O$, there exists a path $\gamma\subset\overline O$ connecting $z$ and $z+\frac aN$ whose length is proportional to $1/N$. Hence,
\begin{equation}
\label{est-int}
\map\left(z+\frac aN\right)-\map(z)-\frac aN\map^\prime(z) = \int_\gamma\big(\map^\prime(t)-\map^\prime(z)\big)dt = o\left(\frac1N\right),
\end{equation}
where the estimate holds uniformly for $z\in T$ and locally uniformly
for $a,b\in\C$. As $|\map(z)|=1$. This means that  
\begin{equation}
\label{Phiab}
\map\left(z+\frac aN\right)\overline{\map\left(z+\frac bN\right)} = 1 + \frac{\tau(a,z)+\overline{\tau(b,z)}}{N} + o\left(\frac1N\right)
\end{equation}
and 
\begin{equation}
\label{Phinab}
\left[\map\left(z+\frac aN\right)\overline{\map\left(z+\frac bN\right)}\right]^N = \exp\left\{\tau(a,z)+\overline{\tau(b,z)}+o(1)\right\},
\end{equation}
where $\tau(\cdot,\cdot)$ was defined before \eqref{omegaaz} and
$o(1)$ is again uniform for $z\in T$ and $a,b$ in compact subsets of
$\C$. As before, we can assume without loss of generality that
$\tau(a,z)+\overline{\tau(b,z)}\neq0$. Then we get from
\eqref{kernel2}, \eqref{Phiab}, \eqref{Phinab}, and the continuity of
$\map^\prime$ that  
\begin{align}
\label{Limit2}
& \lim_{N,s\to\infty}K_{N,s}\left(z+\frac aN,z+\frac bN\right)N^{-2} = \\
& \hspace{3cm}
\frac{|\map^\prime(z)|^2}{\pi}\left(\frac{1-\ell}{2}H_0\left(\tau(a,z)+\overline{\tau(b,z)}\right) 
  + \frac\ell6 H_1\left(\tau(a,z)+\overline{\tau(b,z)}\right)\right) \nonumber
\end{align}
uniformly for $z\in T$ and $a,b$ in compact subsets of $\C$. The limit in \eqref{kernel4} now follows from \eqref{Limit1}, \eqref{Limit2}, and \eqref{H_ell}. 

As obvious from \eqref{tildekernel}, to prove \eqref{kernel5}, it suffices to show that 
\begin{equation}
\label{Limit3}
\lim_{N,s\to\infty}P_K^{-s}(z + a/N)  = \omega(a,z)
\end{equation}
uniformly for $z\in T$ and locally uniformly for $a\in\C$. To this end
observe that an outward normal to $T$ at $z$ is given by
$\map(z)/\map^\prime(z)$.  Hence, the angle between the vectors $a/N$
and $\map(z)/\map^\prime(z)$ is less than $\pi/2$ if and only if the
vector $a\map^\prime(z)/\map(z)=\tau(a,z)$ belongs to the right
half-plane. That is, if $\re(\tau(a,z))>0$. Hence, the limit in
\eqref{Limit3} holds for $\re(\tau(a,z))<0$ as $P_K(z+a/N)\equiv1$ for
such $a$. Moreover, when $\re(\tau(a,z))>0$, \eqref{Limit3} follows
immediately from \eqref{Phinab}. The case $\re(\tau(a,z))=0$ can be
deduced by continuity and the uniformity of the estimate follows from
the uniform character of the estimate in \eqref{est-int}. Finally, the
same arguments yield \eqref{kernel5} for $\ell=0$. \qed 

\appendix

\section{Plots of Correlation Functions}

To provide intuition for the results reported here we consider the
scaled limit of $R_1$ and $R_2$ of the entropic (potential theoretic,
with $K = \overline\D$) ensemble in a 
neighborhood of a point on the unit circle.  By the radial symmetry of the weight it
suffices to restrict ourselves to a neighborhood of 1.  In this case 
$\tau(a,1) = a$ and $\omega(a,1) = \displaystyle e^{-
  \re(a)/\ell}$ for $\re(a)>0$ and $\omega(a,1)=1$ otherwise. As
before, if $\ell=0$ and $\re(a) > 0$ then we take $\omega(a,1) = 0$.
For convenience we define the scaled kernel at $z=1$ by   
\[
\widetilde{H}_{\ell}(a, b) = \omega(a,1) \omega(b,1)  H_\ell\left(a + \overline b\right).
\]
The limiting density of scaled eigenvalues is then given by 
\[
R_1^{\ell}(a) = \widetilde H_{\ell}(a,a),
\]
and the scaling limit of the second correlation function is given by
\[
R_2^{\ell}(a,b) = \widetilde H_{\ell}(a,a) \widetilde H_{\ell}(b,b) - \widetilde H_{\ell}(a,b) \widetilde H_{\ell}(b,a)
\]
The visualizations provided here are for the cases where $a$ and $b$
are either real or on the imaginary axis.  

\subsection{Tangent to the curve}

The tangent line of
the circle at $z=1$ is parallel to the imaginary axis and the local
density of eigenvalues in this direction is given by
\[
\widetilde{H}_{\ell}(i t, i t) = 1.
\]
This is expected since the spatial density of eigenvalues on the unit
circle must be invariant under rotation (and locally, this rotation is given by
translation up the imaginary axis).  

\begin{figure}[h!]
  \centering
  \includegraphics[scale=.5]{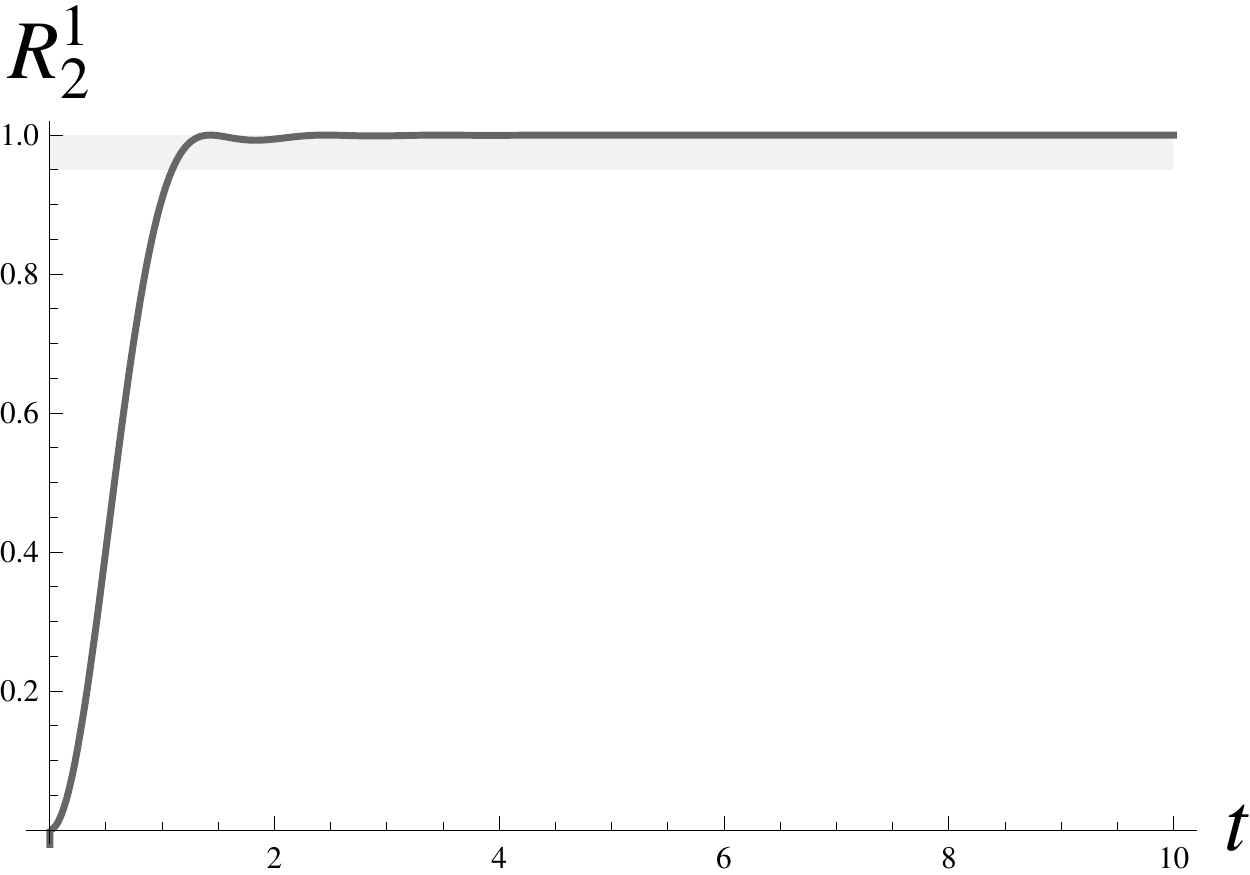}
  \includegraphics[scale=.5]{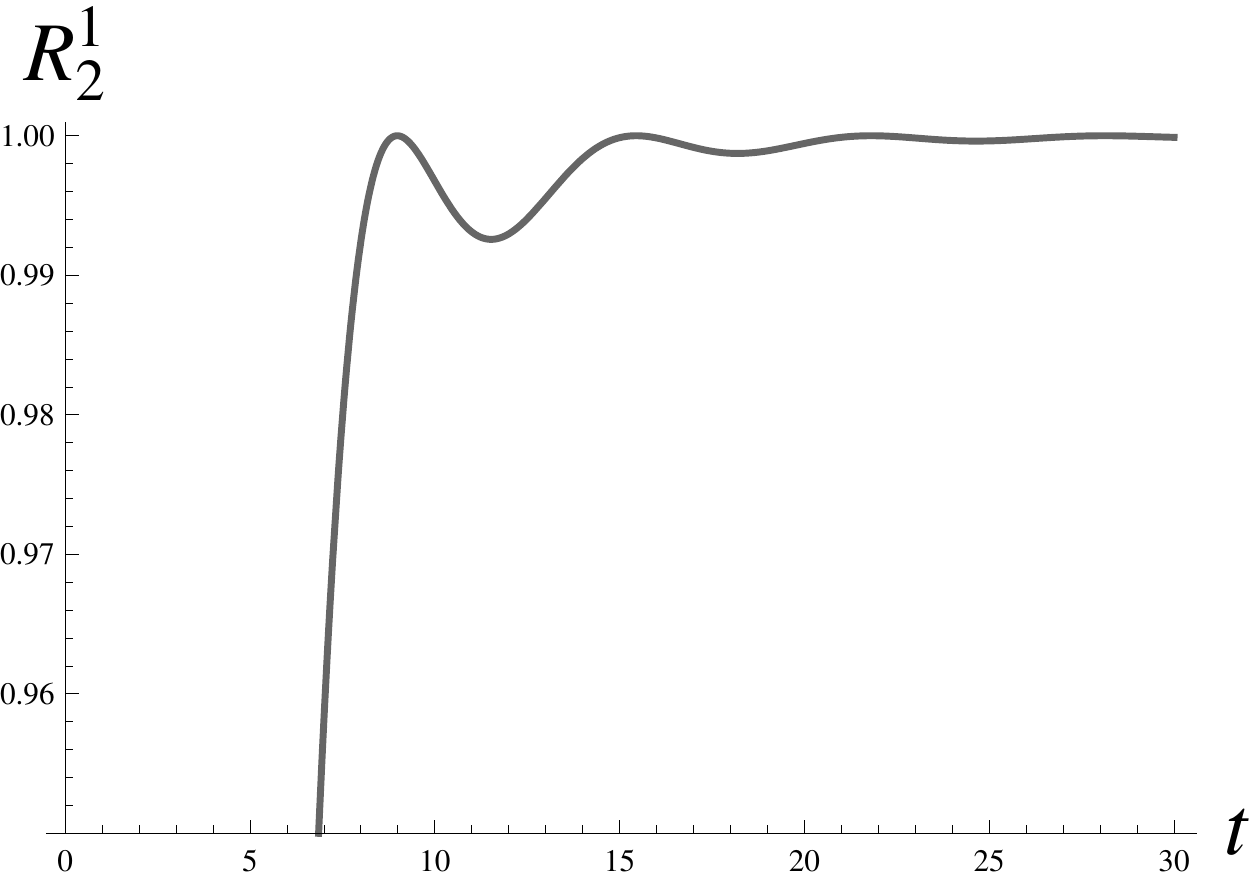}
  \caption{plots of $R^1_2(a,b)$ as a function of $t = -i (a +
    \overline b)$. The second plot is an enlargement of the shaded region.}
\label{fig:1}
\end{figure}

\begin{figure}[h!]
  \centering
  \includegraphics[scale=.5]{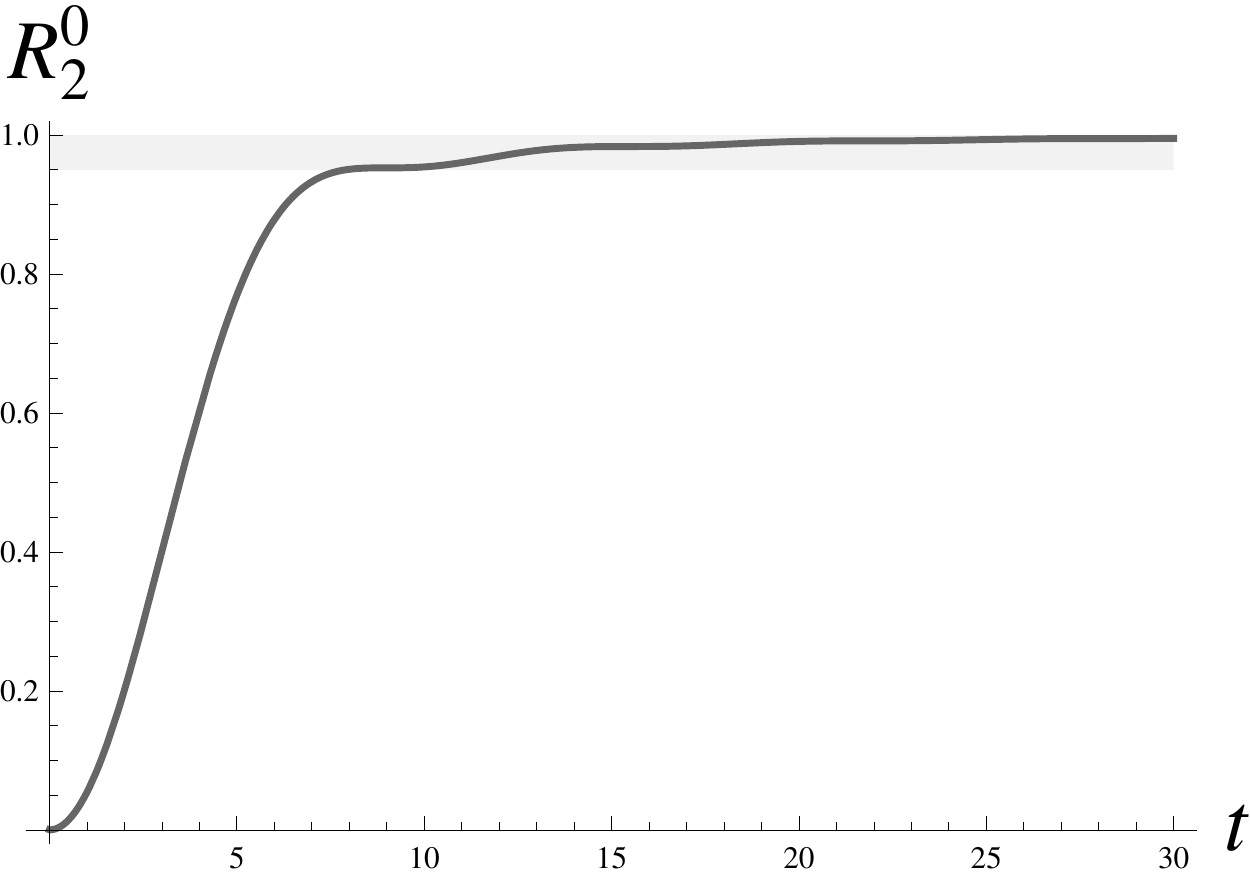}
  \includegraphics[scale=.5]{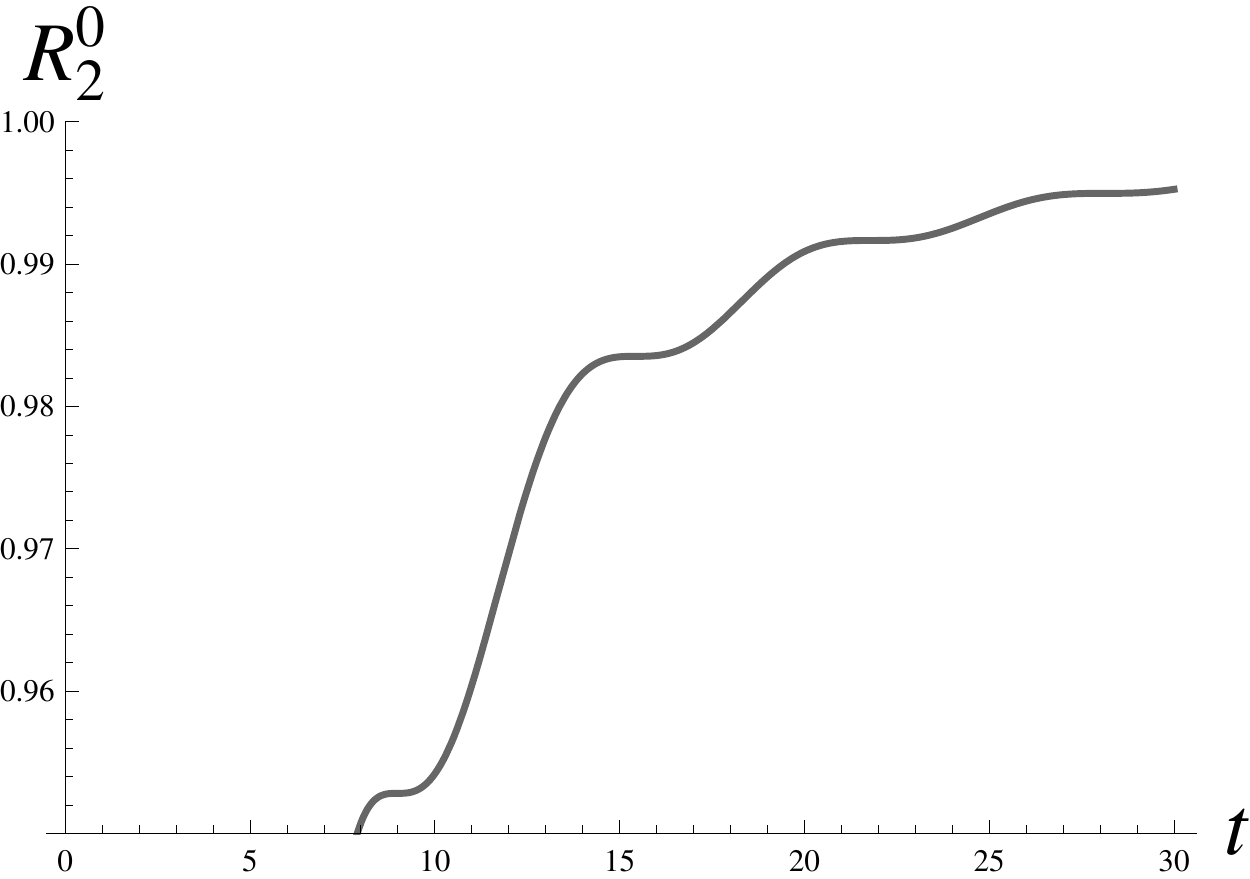}
  \caption{plots of $R^0_2(a,b)$ as a function of $t = -i (a +
    \overline b)$.}
\label{fig:2}
\end{figure}

\begin{figure}[h!]
  \centering
  \includegraphics[scale=1]{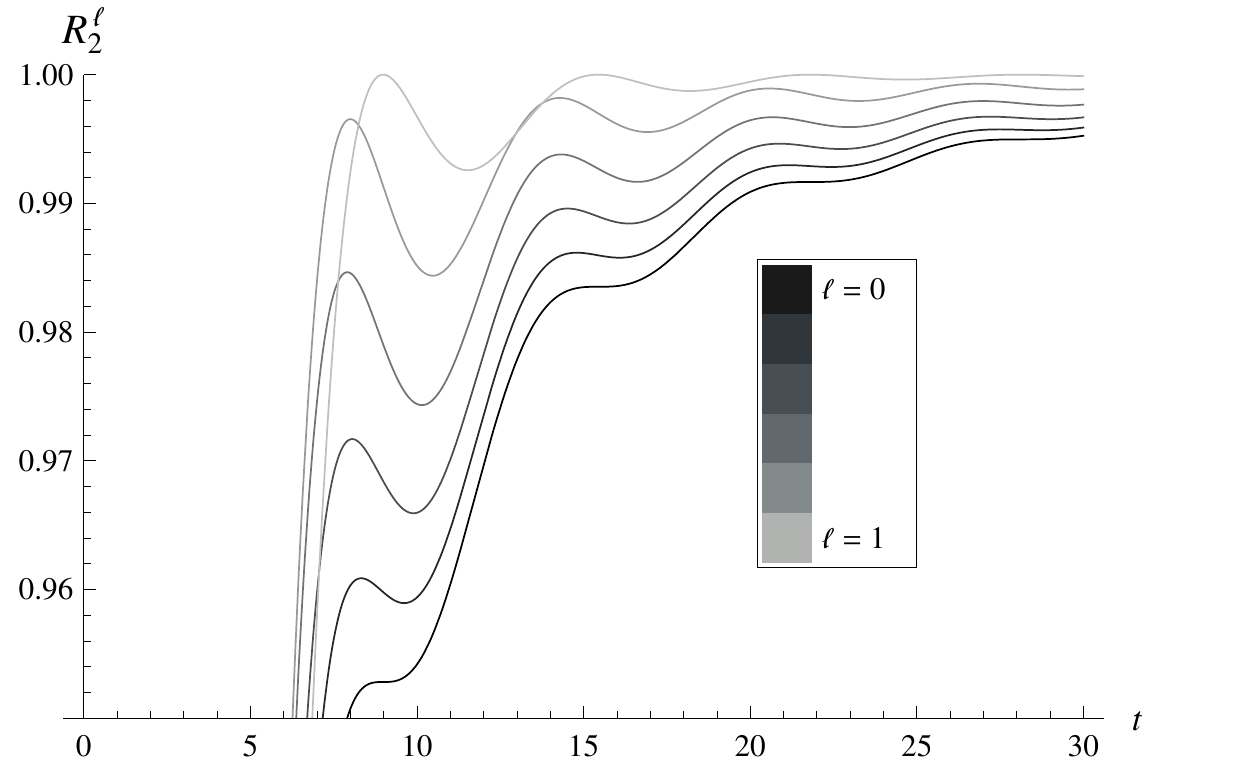}
  \caption{Plot of the interpolation between $R^{0}_2(a,b)$ and
    $R^{1}_2(a,b)$ as a function of $t = -i(a + \overline b)$.} 
\label{fig:3}
\end{figure}

Looking at the second correlation function, when $a$ and $b$ are on
the imaginary axis, we see that $\widetilde{H}_{\ell}$ is a function
of  $t = -i(a + \overline b)$.   Figures~\ref{fig:1}, \ref{fig:2} and \ref{fig:3} show plots of
the second correlation function for various values of $\ell$ in
various regions as a function of $t$.  

By way of comparison we also provide plots of the second scaled
correlation function for ensembles with the sine kernel.
Specifically, 
\[
S(a,b) = 2 \frac{\sin\big((a - b)/2\big)}{(a - b)} \qquad \mbox{and} \qquad
R_2^{\sin}(a,b) = 1 - S(a,b)^2.
\]
(The slightly unusual normalization given by the superfluous appearing
factors of 2 in the first equation arises in the scaling limit when we take the
expected distance between eigenvalues to be $2 \pi$---this allows for
the most accurate comparison with our other figures).  

\begin{figure}[h!]
  \centering
  \includegraphics[scale=.5]{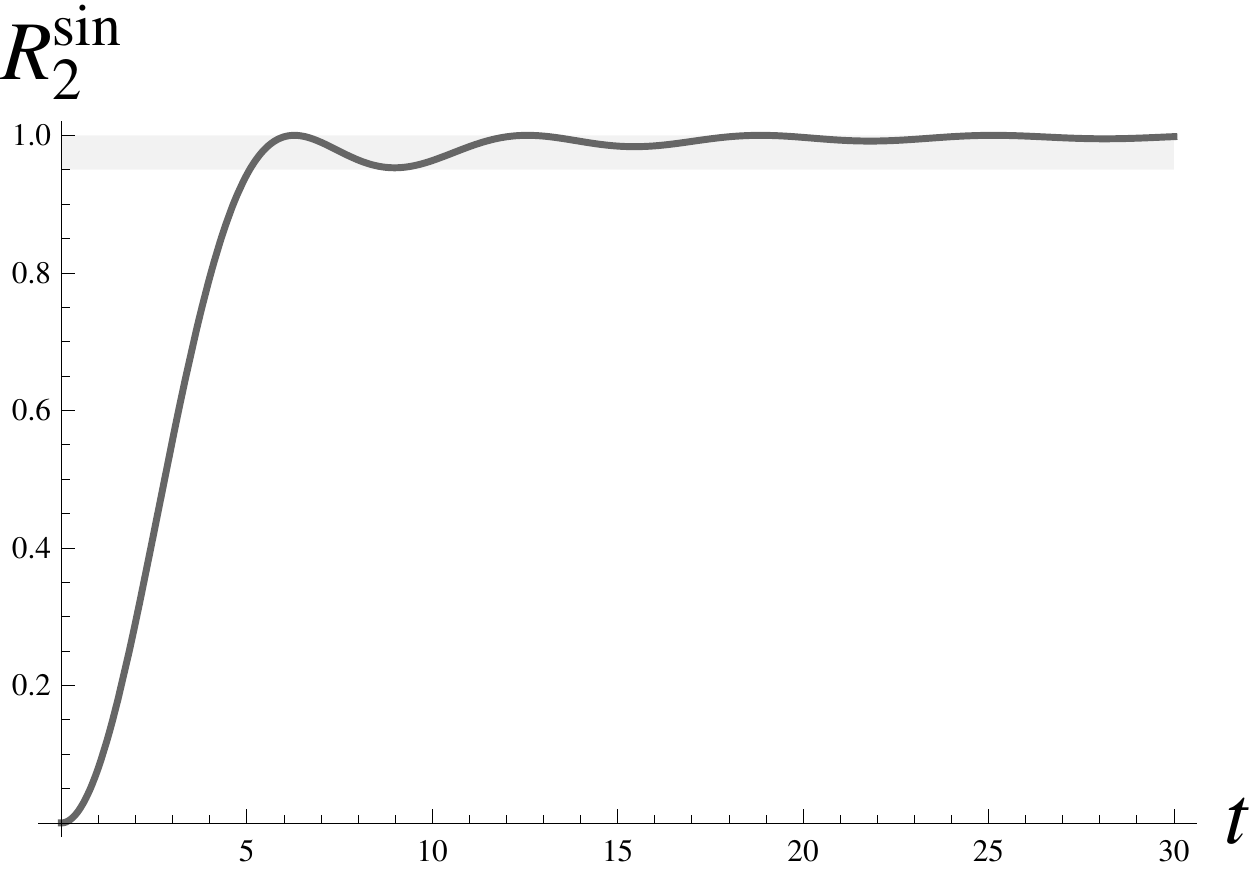}
  \includegraphics[scale=.5]{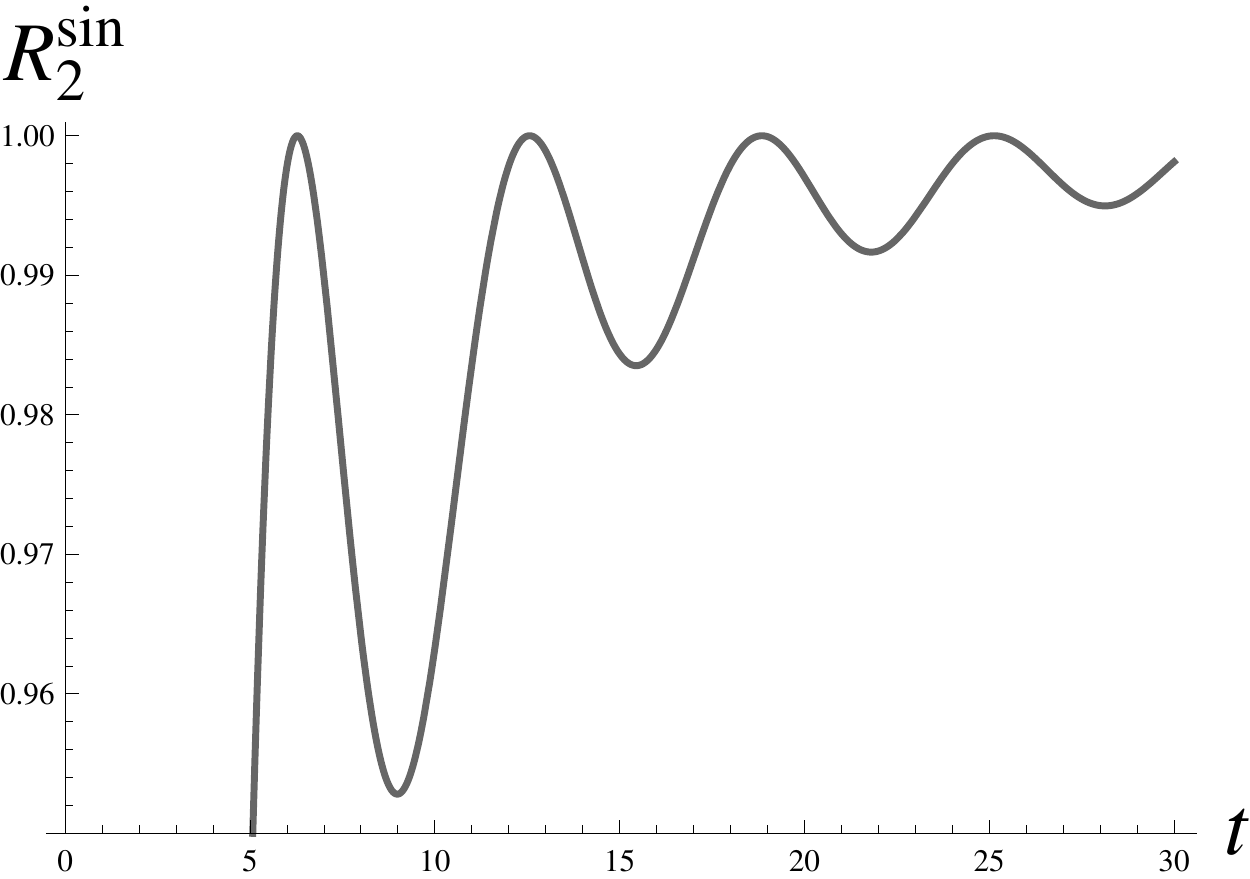}
  \caption{plots of $R^{\sin}_2(a,b)$ as a function of $t = a - b$.}
\label{fig:4}
\end{figure}

\subsection{Normal to the curve}

In the regime where $a$ and $b$ are real, we are looking in a
neighborhood of $z = 1$ in a direction perpendicular to that where the
density of eigenvalues becomes constant.  That is, the first scaled
correlation function should decay as $a$ moves away from 0.  Negative
$a$ corresponds to moving into $K$ (where the potential is
constant) whereas positive $a$ corresponds to moving away from $K$
where the potential acts to make $\infty$ 
repulsive.  As $\ell$ decreases to 0, the field increases in strength
until at $\ell = 0$ there is no possibility that an eigenvalue can be
outside $K$.  That is, when $\ell = 0$ the first
correlation function vanishes for $a > 0$.  

\begin{figure}[h!]
  \centering
  \includegraphics[scale=1]{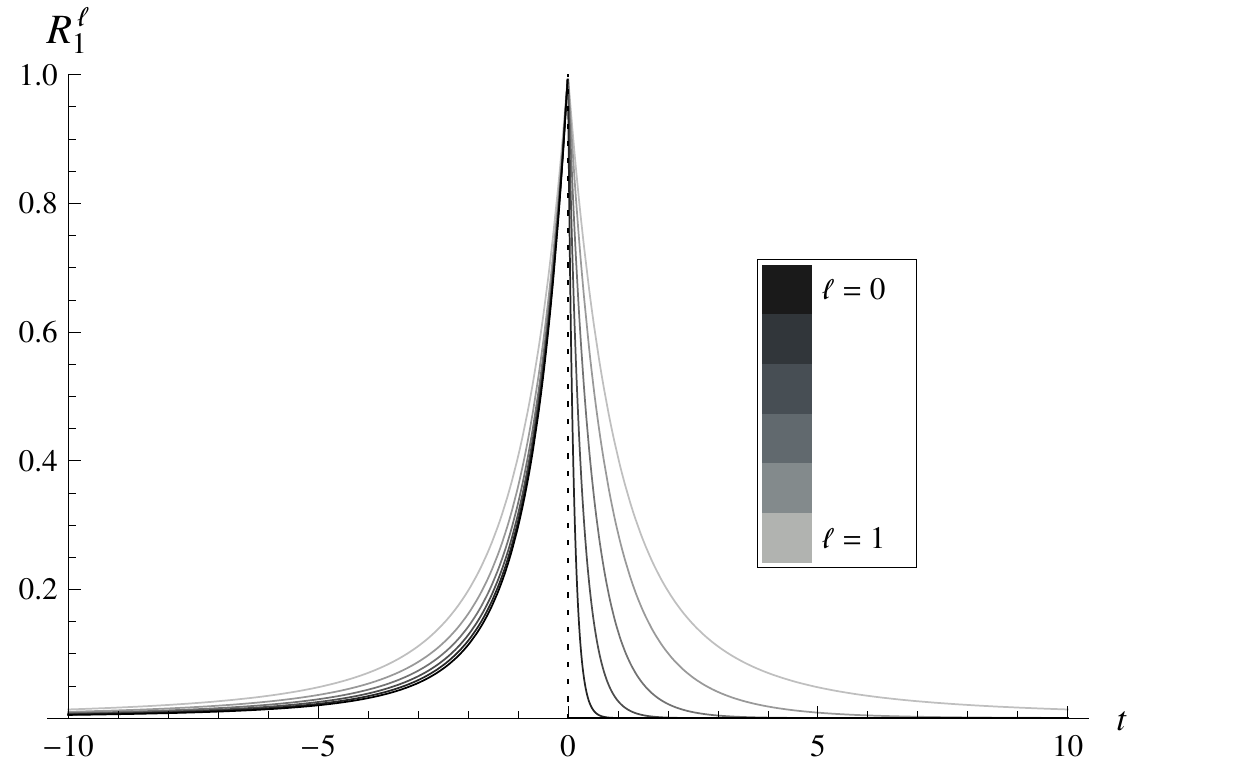}
  \caption{Plot of the interpolation between $R^{0}_1(a)$ and
    $R^{1}_1(a)$ as a function of $t = -\mathrm{Re}(a)$.  When
    $\ell=0$ there is a sharp cutoff at $t = 0$.  } 
\label{fig:5}
\end{figure}

When $a$ and $b$ are real, $R_2^{\ell}(a, b)$ is no longer a function
of a linear combination of $a$ and $b$, and we plot this as a surface
for $\ell = 0$ and $\ell = 1$.  

\begin{figure}[h!]
  \centering
  \includegraphics[scale=.3]{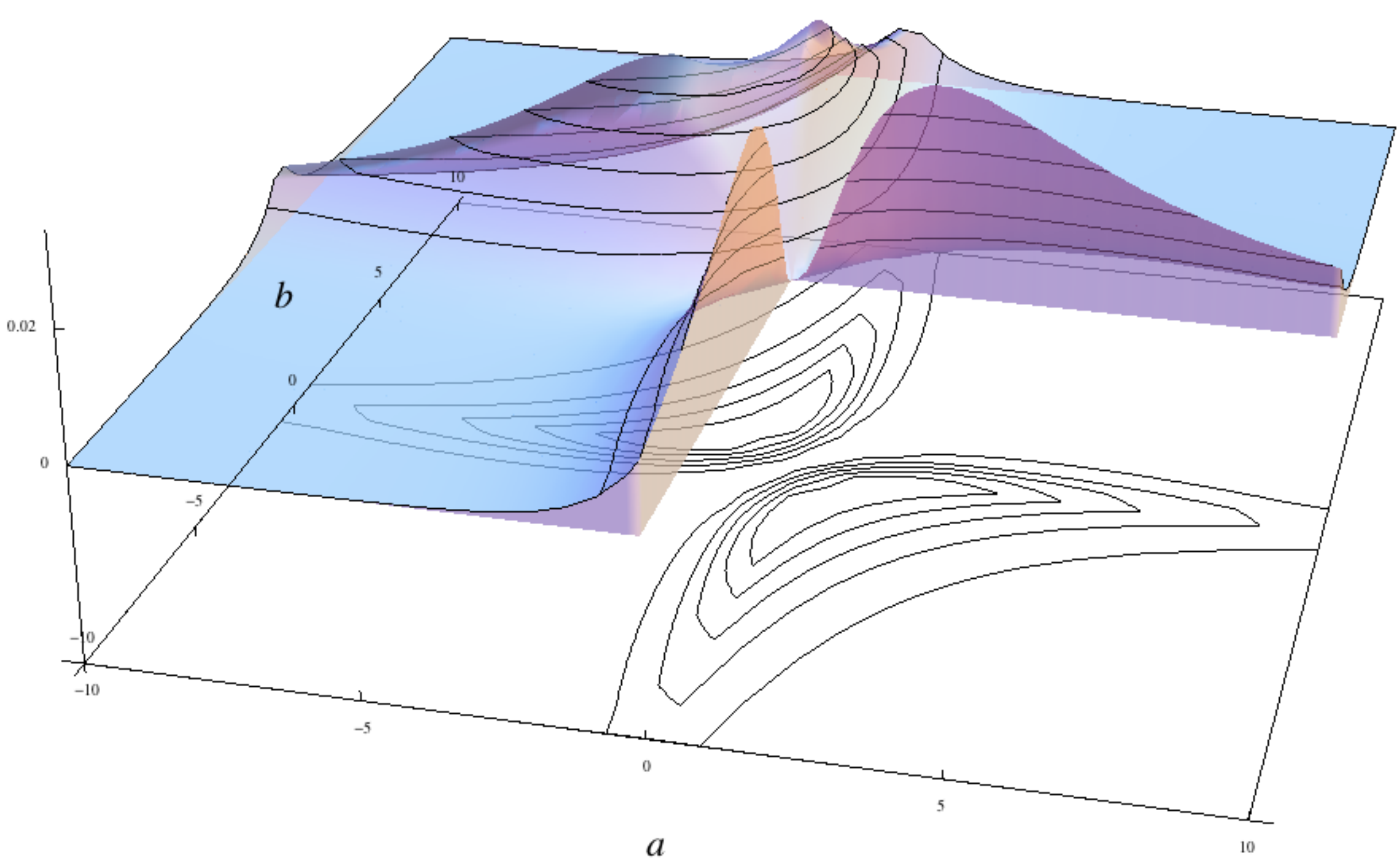}
  \caption{Plot of $R^{1}_2(a,b)$ when $a$ and $b$ are real, with part
  of the surface removed to see the cross-section.} 
\label{fig:6}
\end{figure}

\begin{figure}[h!]
  \centering
  \includegraphics[scale=.8]{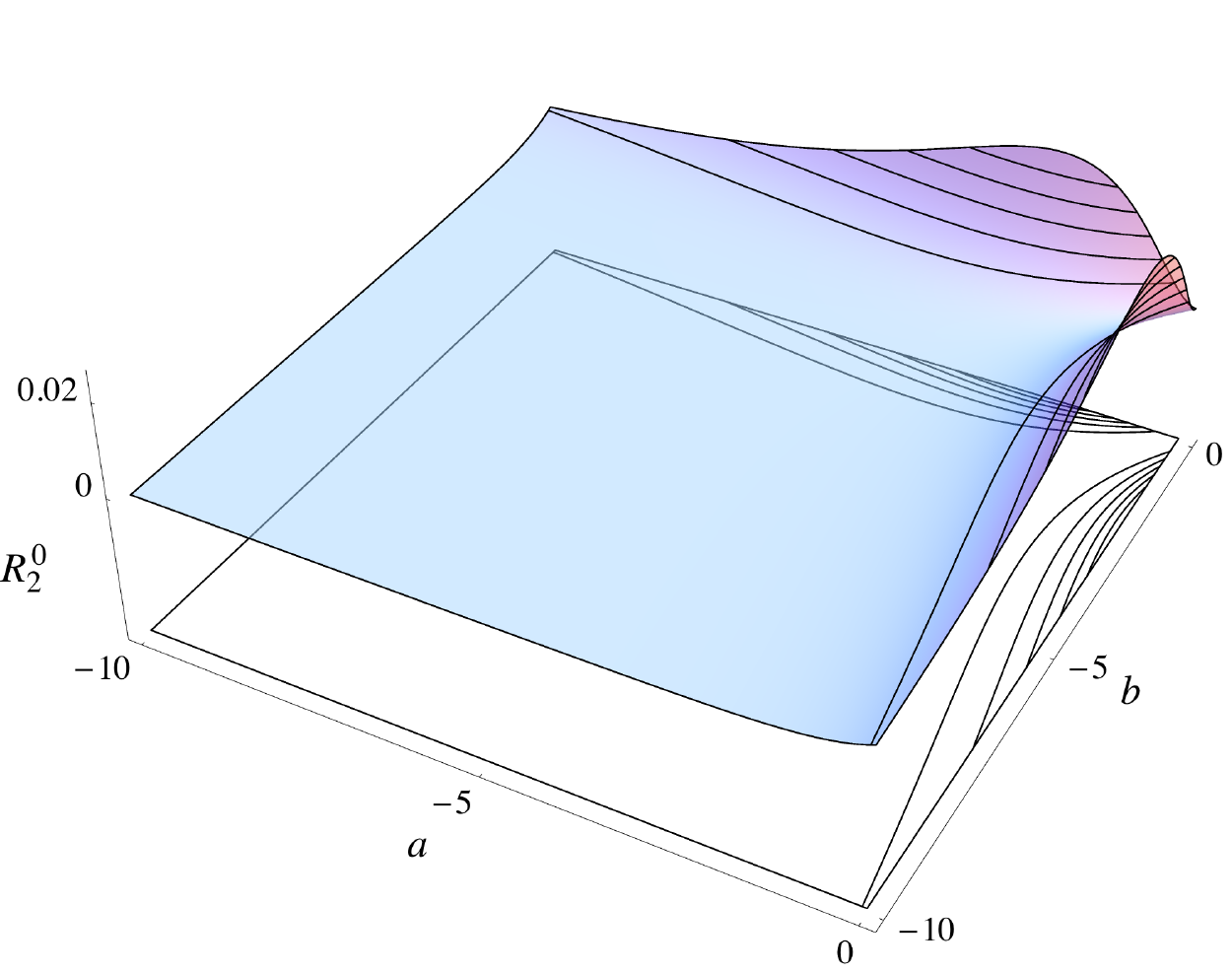}
  \caption{Plot of $R^{0}_2(a,b)$ when $a$ and $b$ are real. Note that
    $R^{0}(a,b)$ is identically zero if either $a$ or $b$ is greater
    than 0.} 
\label{fig:7}
\end{figure}


\begin{thebibliography}{10}

\bibitem{MR2760897}
G.W. Anderson, A. Guionnet, and O. Zeitouni.
\newblock {\em An introduction to random matrices}, volume 118 of {\em
  Cambridge Studies in Advanced Mathematics}.
\newblock Cambridge University Press, Cambridge, 2010.

\bibitem{Carl23}
T.~Carleman.
\newblock {\"U}ber die {A}pproximation analytisher {F}unktionen durch lineare
  {A}ggregate von vorgegebenen {P}otenzen.
\newblock {\em Ark. Mat. Astr. Fys.}, 17(9):1--30, 1922.

\bibitem{LingLieChau|YueYu1992452}
L.-L. Chau and Y. Yu.
\newblock Unitary polynomials in normal matrix models and wave functions for
  the fractional quantum hall effects.
\newblock {\em Physics Letters A}, 167(5-6):452 -- 458, 1992.

\bibitem{MR1643533}
L.-L. Chau and O. Zaboronsky.
\newblock On the structure of correlation functions in the normal matrix model.
\newblock {\em Comm. Math. Phys.}, 196(1):203--247, 1998.

\bibitem{MR1868596}
S.-J. Chern and J.D. Vaaler.
\newblock The distribution of values of {M}ahler's measure.
\newblock {\em J. Reine Angew. Math.}, 540:1--47, 2001.

\bibitem{DragMina10}
P.~Dragnev and E.~Mi{\~n}a-D\'iaz.
\newblock Asymptotic behavior and zero distribution of {C}arleman orthogonal
  polynomials.
\newblock {\em J. Approx. Theory}, 162(11):1982--2003, 2010.

\bibitem{DragMina10a}
P.~Dragnev and E.~Mi{\~n}a-D\'iaz.
\newblock On series representation for {Carleman} orthogonal polynomials.
\newblock {\em Proc. Amer. Math. Soc.}, 138(12):4271--4279, 2010.

\bibitem{Gaier}
D.~Gaier.
\newblock {\em Lectures on Complex Approximation}.
\newblock Birkha\"user Boston, Inc., 1987.

\bibitem{Gakhov}
F.D. Gakhov.
\newblock {\em Boundary Value Problems}.
\newblock Dover Publications, Inc., New York, 1990.

\bibitem{MR0173726}
J. Ginibre.
\newblock Statistical ensembles of complex, quaternion, and real matrices.
\newblock {\em J. Mathematical Phys.}, 6:440--449, 1965.

\bibitem{GohbergGoldbergKrupnik}
I.~Gohberg, S.~Goldberg, and N.~Krupnik.
\newblock {\em Traces and Determinants of Linear Operators}, volume 116 of {\em
  Operator Theory Advances and Applications}.
\newblock Birkh\"auser, Boston, 2000.

\bibitem{Hormander}
L.~H\"ormander.
\newblock {\em An Introduction to Complex Analysis in Several Variables}.
\newblock The University Series in Higher Mathematics. D. Van Nostrand Company,
  Inc., Princeton, New Jersey, 1966.

\bibitem{Lub10}
D.S. Lubinsky.
\newblock Universality type limits for {Bergman} orthogonal polynomials.
\newblock {\em Comput. Methods Funct. Theory}, 10:135--154, 2010.

\bibitem{Mehta1960420}
M.L. Mehta and M.~Gaudin.
\newblock On the density of eigenvalues of a random matrix.
\newblock {\em Nuclear Physics}, 18:420 -- 427, 1960.

\bibitem{Mina08}
E.~Mi{\~n}a-D\'iaz.
\newblock An asymptotic integral representation for {C}arleman orthogonal
  polynomials.
\newblock {\em Int. Math. Res. Not. IMRN}, (16), 2008.
\newblock Art. ID rnn065, 38pp.

\bibitem{PhysRevE.55.205}
G. Oas.
\newblock Universal cubic eigenvalue repulsion for random normal matrices.
\newblock {\em Phys. Rev. E}, 55(1):205--211, Jan 1997.

\bibitem{Ransford}
T.~Ransford.
\newblock {\em Potential Theory in the Complex Plane}, volume~28 of {\em London
  Mathematical Society Student Texts}.
\newblock Cambridge University Press, Cambridge, 1995.

\bibitem{MR2145532}
C.D. Sinclair.
\newblock The distribution of {M}ahler's measures of reciprocal polynomials.
\newblock {\em Int. J. Math. Math. Sci.}, (49-52):2773--2786, 2004.

\bibitem{MR2407817}
C.D. Sinclair.
\newblock The range of multiplicative functions on {$\Bbb C[x],\ \Bbb R[x]$}
  and {$\Bbb Z[x]$}.
\newblock {\em Proc. Lond. Math. Soc. (3)}, 96(3):697--737, 2008.

\bibitem{MR2707617}
C.D. Sinclair.
\newblock {\em Multiplicative distance functions}.
\newblock ProQuest LLC, Ann Arbor, MI, 2005.
\newblock Thesis (Ph.D.)--The University of Texas at Austin.

\bibitem{Suetin}
P.K. Suetin.
\newblock {\em Polynomials Orthogonal over a Region and Bieberbach
  Polynomials}, volume 100 of {\em Proc. Steklov Inst. Math.}
\newblock Amer. Math. Soc. Translations, 1974.

\bibitem{TW}
C.A. Tracy and H. Widom.
\newblock Correlation functions, cluster functions, and spacing distributions
  for random matrices.
\newblock {\em J. Statist. Phys}, 92:809--835, 1998.

\bibitem{springerlink:10.1007/BF01040581}
S.A. Yuzvinskii.
\newblock Computing the entropy of a group of endomorphisms.
\newblock {\em Siberian Mathematical Journal}, 8:172--178, 1967.
\newblock 10.1007/BF01040581.



\end{thebibliography}
\end{document}